\documentclass[a4paper,12pt,oneside,reqno]{amsart}
\usepackage{comment,amsmath,amssymb,amsthm,pstricks,changepage,pifont,fullpage,todonotes,hyperref,pstricks-add}
\usepackage[english]{babel}

\newtheorem{theorem}{Theorem}[section]

\newtheorem{lemma}[theorem]{Lemma}
\newtheorem{corollary}[theorem]{Corollary}

\newtheorem{definition}[theorem]{Definition}

\DeclareMathOperator{\lsize}{{ls}}
\DeclareMathOperator{\conv}{conv}

\begin{document}

\title{Linear pencils encoded in the Newton polygon}
\author{Wouter Castryck and Filip Cools}

\maketitle

\begin{abstract}
  \noindent Let $C$ be an algebraic curve defined by a sufficiently generic bivariate Laurent polynomial with given Newton polygon $\Delta$.
   It is classical that the geometric genus of $C$ equals the number of lattice points in the interior of $\Delta$.
   In this paper we give similar combinatorial interpretations for the gonality, the Clifford index and the Clifford dimension, 
   by removing a technical assumption from a recent result of Kawaguchi.
   More generally, the method shows that apart from certain well-understood exceptions, every base-point free
   pencil whose degree equals or slightly exceeds the gonality is combinatorial, in the sense
   that it corresponds to projecting $C$ along a lattice direction. 
   Along the way we prove various features of combinatorial pencils. For instance, we give an interpretation
   for the scrollar invariants associated to a combinatorial pencil, and show how one can tell whether
   the pencil is complete or not. 
   
   Among the applications, we find that every smooth projective curve admits at most
   one Weierstrass semi-group of embedding dimension $2$, and that if a non-hyperelliptic
   smooth projective curve $C$ of genus $g \geq 2$ can be embedded in the $n^\text{th}$ Hirzebruch surface $\mathcal{H}_n$,
   then $n$ is actually an invariant of $C$.\\
   
\end{abstract}

\setcounter{tocdepth}{1}
\tableofcontents
\vspace{-1.13cm}
Accompanying Magma files\footnote{Available at \url{http://users.ugent.be/~wcastryc/}}: \verb"basic_commands.m", \verb"gonal.m", \verb"neargonal.m"

\newpage

\section{Introduction}

Let $k$ be an algebraically closed field of characteristic zero, let $\mathbb{T}^2 = (k^\ast)^2$ be the $2$-dimensional
torus over $k$, and let $f \in k[x^{\pm 1}, y^{\pm 1} ]$ be
an irreducible Laurent polynomial. Denote by $U(f)$ the curve in $\mathbb{T}^2$ defined by $f$.
Let $\Delta(f) \subset \mathbb{R}^2$ be the Newton polygon of $f$, which we
always assume to be two-dimensional. We say that $f$ is \emph{non-degenerate with respect to its Newton polygon}
if for every face $\tau \subset \Delta(f)$ (including $\Delta(f)$ itself) the system
\[ f_\tau = \frac{\partial f_\tau }{\partial x} = \frac{\partial f_\tau}{\partial y} = 0 \]
has no solutions in $\mathbb{T}^2$. (Here $f_\tau$ is obtained from $f$ by only considering
those terms that are supported on $\tau$.) For a two-dimensional lattice polygon $\Delta \subset \mathbb{R}^2$, we say
that $f$ is \emph{$\Delta$-non-degenerate} if it is non-degenerate with respect to its Newton polygon and
$\Delta(f) = \Delta$.
For Laurent polynomials that are supported on $\Delta$,
the condition of $\Delta$-non-degeneracy is generically satisfied, in the sense
that it is characterized by the non-vanishing of
\[ \text{Res}_\Delta \left(f, x\frac{\partial f}{\partial x}, y \frac{\partial f}{\partial y} \right)  \]
(where $\text{Res}_\Delta$ is the sparse resultant; see \cite[Prop.\,1.2]{CastryckVoight} and \cite[Thm.\,10.1.2]{GKZ} for an according discussion).
An algebraic curve $C / k$
is called \emph{$\Delta$-non-degenerate} if it is birationally equivalent to $U(f)$ for some $\Delta$-non-degenerate Laurent
polynomial $f \in k[x^{\pm 1}, y^{\pm 1}]$.\\

\noindent \emph{Remark.}
Sometimes in the existing literature a projectively embedded variety is called \emph{non-degenerate} if it is not contained
in a hyperplane. Our notion of non-degeneracy is unrelated to this.\\

It is well-known that if $C$ is $\Delta$-non-degenerate, then several of its geometric properties
are encoded in the combinatorics of $\Delta$. The most prominent example is that the geometric genus equals the
number of lattice points in the interior of $\Delta$ \cite[\S4~Ass.\,2]{Khovanskii}. The proof of this fact 
is briefly recalled at the beginning of Section~\ref{section_nondegistoric}, because it
entails an explicit description of the canonical map that will play a role in Section~\ref{section_maroni}.
Other known examples are that one can tell from
$\Delta$ whether $C$ is hyperelliptic or not \cite[Lem.\,3.2.9]{Koelman}, and whether it is trigonal or not
\cite[Lem.\,3]{CaCo}. Recently, this was extended to arbitrary gonalities by Kawaguchi \cite[Thm.\,1.3]{Kawaguchi}
under the technical assumption that $C$ is not birationally equivalent to a smooth plane projective curve.
 
In Section~\ref{section_gonality} we 
revisit Kawaguchi's proof, while making a more explicit
connection with the language of Newton polygons and getting rid of the above technical assumption. 
Kawaguchi's method yields that apart from some well-understood exceptional instances
of $\Delta$, every gonality pencil on $C$
is \emph{combinatorial}, in the sense that it corresponds to a projection of the form $(x,y) \mapsto x^ay^b$ for coprime $a,b \in \mathbb{Z}$.
In this case, the gonality is easily seen to
equal the \emph{lattice width} of $\Delta$ (this notion will be recalled in Section~\ref{section_combpencils}).
This settles a conjecture by the current authors \cite[Conj.\,1]{CaCo}, although most cases, including
all lattice polygons whose number of interior lattice points is not of the form $(d-1)(d-2)/2$, were
already covered by Kawaguchi's work.

In Section~\ref{section_neargonal} we apply the same method to near-gonal pencils, i.e.\ 
base-point free linear systems of the form $g^1_{\gamma + 1}$, where $\gamma$ is the gonality. 
It again turns out that, apart from some reasonably well understood exceptions, every such pencil
is combinatorial.

Then in Section~\ref{section_clifford}, we prove that also the Clifford index and the Clifford dimension of $C$ are fully
determined by the combinatorics of $\Delta$. This is again inspired by \cite{Kawaguchi}, but thanks to
our coverage of the case of smooth projective plane curves (i.e., curves of Clifford dimension $2$)
we are able to fill in the missing spots.
In particular, we obtain a purely combinatorial criterion for determining whether $C$ is birationally equivalent to
a smooth projective curve in $\mathbb{P}^2$ or not. 

Note that, as an immediate corollary to all this, we obtain that the gonality, the Clifford index and the Clifford dimension do not 
depend on the specific choice of our $\Delta$-non-degenerate curve $C$. This is an extension to arbitrary toric surfaces of a
recent theorem by Lelli-Chiesa \cite[Thm.\,1.2]{LelliChiesa} on families of curves on rational (e.g.\ toric) surfaces that carry an anticanonical pencil. 


Next, 
in Section~\ref{section_maroni}, we show that the scrollar invariants associated
to a combinatorial pencil (which specialize to the classical Maroni invariants in the case of a $g^1_3$) have a natural 
combinatorial interpretation. The same interpretation allows one to decide whether a given combinatorial pencil is complete or not.


Finally, Section~\ref{section_applications} discusses a number of applications.
One potential use of our results is as a tool for constructing examples of curves having
certain prescribed invariants (and for finding
lower bounds on the dimension of the corresponding moduli space).
Among the other byproducts we find that 
\begin{itemize}
\item any curve (not necessarily non-degenerate) admits at most one Weierstrass semi-group of embedding dimension two,
\item if $C$ is a non-hyperelliptic smooth projective curve of genus $g \geq 2$ in the $n^\text{th}$ Hirzebruch surface $\mathcal{H}_n$,
then $n$ is actually an invariant of $C$.
\end{itemize}

\section{Notation, terminology and conventions} \label{section_notation}

For lattice polygons $\Delta, \Delta' \subset \mathbb{R}^2$, we say that $\Delta$ is \emph{equivalent} to $\Delta'$ (notation: $\Delta \cong \Delta'$) if $\Delta'$ is obtained from $\Delta$ through a unimodular transformation, i.e.\
through a transformation of the form
\[ \upsilon: \mathbb{R}^2 \rightarrow \mathbb{R}^2 : \begin{pmatrix} i \\ j \\ \end{pmatrix} \mapsto A \begin{pmatrix} i \\ j \\ \end{pmatrix} +
\begin{pmatrix} a_1 \\ a_2 \\ \end{pmatrix}, \qquad A \in \text{GL}_2(\mathbb{Z}), \ a_1,a_2 \in \mathbb{Z}. \]
If $A$ can be taken the unit matrix, we sometimes write $\Delta \cong_t \Delta'$ to emphasize
that $\Delta$ is obtained from $\Delta'$ through a translation.
Note that if a Laurent polynomial
\[ f = \sum_{(i,j) \in \Delta \cap \mathbb{Z}^2} c_{i,j}(x,y)^{(i,j)} \]
is $\Delta$-non-degenerate (where $(x,y)^{(i,j)}$ means $x^iy^j$) and $\upsilon$ is a unimodular transformation, then
\[ f^\upsilon = \sum_{(i,j) \in \Delta \cap \mathbb{Z}^2} c_{i,j}(x,y)^{\upsilon(i,j)} \]
is $\upsilon(\Delta)$-non-degenerate, and $U(f) \cong U(f^\upsilon)$. (Every unimodular transformation
induces an automorphism of $\mathbb{T}^2$.) 

It is convenient to introduce a special notation for certain recurring polygons:
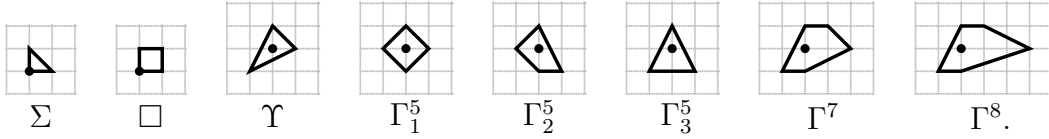
\begin{center}
\psset{unit=0.3cm}
\begin{pspicture}(0,-1.3)(3,3.3)
\psgrid[subgriddiv=1,griddots=10,gridcolor=lightgray,gridlabels=0pt](0,0)(3,3)
\pspolygon[linewidth=1.3pt](1,1)(2,1)(1,2)
\pscircle[fillstyle=solid,fillcolor=black](1,1){0.2}
\rput(1.5,-1){$\Sigma$}
\end{pspicture}
\quad
\psset{unit=0.3cm}
\begin{pspicture}(0,-1.3)(3,3.3)
\psgrid[subgriddiv=1,griddots=10,gridcolor=lightgray,gridlabels=0pt](0,0)(3,3)
\pspolygon[linewidth=1.3pt](1,1)(2,1)(2,2)(1,2)
\pscircle[fillstyle=solid,fillcolor=black](1,1){0.2}
\rput(1.5,-1){$\square$}
\end{pspicture}
\quad
\begin{pspicture}(0,-1.3)(4,4.3)
\psgrid[subgriddiv=1,griddots=10,gridcolor=lightgray,gridlabels=0pt](0,0)(4,4)
\pspolygon[linewidth=1.3pt](3,2)(2,3)(1,1)
\pscircle[fillstyle=solid,fillcolor=black](2,2){0.2}
\rput(2,-1){$\Upsilon$}
\end{pspicture}
\quad
\begin{pspicture}(-1,-2.3)(3,3.3)
\psgrid[subgriddiv=1,griddots=10,gridcolor=lightgray,gridlabels=0pt](-1,-1)(3,3)
\pspolygon[linewidth=1.3pt](0,1)(1,0)(2,1)(1,2)
\pscircle[fillstyle=solid,fillcolor=black](1,1){0.2}
\rput(1,-2){$\Gamma^5_1$}
\end{pspicture}
\quad
\begin{pspicture}(-1,-2.3)(3,3.3)
\psgrid[subgriddiv=1,griddots=10,gridcolor=lightgray,gridlabels=0pt](-1,-1)(3,3)
\pspolygon[linewidth=1.3pt](0,1)(1,0)(2,0)(1,2)
\pscircle[fillstyle=solid,fillcolor=black](1,1){0.2}
\rput(1,-2){$\Gamma^5_2$}
\end{pspicture}
\quad
\begin{pspicture}(-1,-2.3)(3,3.3)
\psgrid[subgriddiv=1,griddots=10,gridcolor=lightgray,gridlabels=0pt](-1,-1)(3,3)
\pspolygon[linewidth=1.3pt](0,0)(2,0)(1,2)
\pscircle[fillstyle=solid,fillcolor=black](1,1){0.2}
\rput(1.2,-2){$\Gamma^5_3$}
\end{pspicture}
\quad
\begin{pspicture}(-1,-2.3)(4,3.3)
\psgrid[subgriddiv=1,griddots=10,gridcolor=lightgray,gridlabels=0pt](-1,-1)(4,3)
\pspolygon[linewidth=1.3pt](0,0)(1,0)(3,1)(2,2)(1,2)
\pscircle[fillstyle=solid,fillcolor=black](1,1){0.2}
\rput(1.8,-2){$\Gamma^7$}
\end{pspicture}
\quad
\begin{pspicture}(-1,-2.3)(5,3.3)
\psgrid[subgriddiv=1,griddots=10,gridcolor=lightgray,gridlabels=0pt](-1,-1)(5,3)
\pspolygon[linewidth=1.3pt](0,0)(1,0)(4,1)(2,2)(1,2)
\pscircle[fillstyle=solid,fillcolor=black](1,1){0.2}
\rput(2.3,-2){$\Gamma^8$.}
\end{pspicture}
\end{center}
Here the bold-marked lattice point indicates the point $(0,0) \in \mathbb{R}^2$, although we are usually interested in lattice polygons up to equivalence only.
Thus $\Sigma$ is the standard simplex, and $d\Sigma$ (Minkowski multiple) is the Newton polygon
of a generic degree $d$ polynomial. If $\Delta$ is a two-dimensional lattice polygon, we denote by $\Delta^{(1)}$ the convex
hull of its interior lattice points. The boundary of $\Delta$ is denoted by $\partial \Delta$.\\

\noindent \emph{Example.} For $d\geq 3$ one has $(d\Sigma)^{(1)} \cong (d-3)\Sigma$. For
$d \geq 1$ one has $(d\Upsilon)^{(1)} \cong (d-1)\Upsilon$.\\

\noindent \emph{Remark.} Occasionally, we will also apply the notation $\Delta^{(1)}$ to convex
polygons $\Delta$ that are lower-dimensional and/or take vertices outside the lattice $\mathbb{Z}^2$. Here again we mean
the convex hull of the lattice points in the interior of $\Delta$, where the interior is understood to be empty in the lower-dimensional case.\\

If $\Delta^{(1)}$ is two-dimensional, then the
set of lattice polygons $\Gamma$ for which $\Gamma^{(1)} = \Delta^{(1)}$ admits
a maximum with respect to inclusion \cite[Lem.\,9]{HaaseSchicho}. We denote this maximum by $\Delta^\text{max}$.
It can be characterized as follows.
Write $\Delta^{(1)}$ as an intersection of half-spaces
\[ \bigcap_{\ell=1}^r H_\ell, \qquad \text{with $H_\ell = \left\{ \, \left. (i,j) \in \mathbb{R}^2 \, \right| \, <\! (i,j), v_\ell \!> \ \geq - a_\ell \, \right\}$,} \]
where 
$< \! \cdot , \cdot \!>$ denotes the standard inner product on $\mathbb{R}^2$ and
$v_1, \dots, v_r$ are primitive inward pointing normal vectors of the edges of $\Delta^{(1)}$.
Then
\[ \Delta^{\text{max}} = \bigcap_{\ell=1}^r H_\ell^{(-1)}, \qquad \text{where $H_\ell^{(-1)} = \left\{ \, \left. (i,j) \in \mathbb{R}^2 \, \right| \, <\! (i,j), v_\ell \!> \ \geq - a_\ell - 1 \, \right\}$.} \]
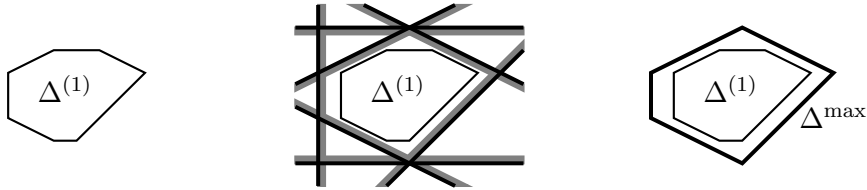
\begin{center}
\psset{unit=0.3cm}
\begin{pspicture}(-1,-1)(9,7)
\pspolygon(3,1)(4,1)(7,4)(5,5)(3,5)(1,4)(1,2)
\rput(3.5,3.3){\small $\Delta^{(1)}$}
\end{pspicture}
\qquad \quad
\begin{pspicture}(-1,-1)(9,7)
\pspolygon(3,1)(4,1)(7,4)(5,5)(3,5)(1,4)(1,2)
\pspolygon[linecolor=gray,fillcolor=gray,fillstyle=solid](0,-1)(0,7)(0.3,7)(0.3,-1)
\pspolygon[linecolor=gray,fillcolor=gray,fillstyle=solid](-1,3.5)(6,7)(6.6,7)(-1,3.2)
\pspolygon[linecolor=gray,fillcolor=gray,fillstyle=solid](-1,6)(9,6)(9,5.7)(-1,5.7)
\pspolygon[linecolor=gray,fillcolor=gray,fillstyle=solid](9,3.5)(2,7)(1.4,7)(9,3.2)
\pspolygon[linecolor=gray,fillcolor=gray,fillstyle=solid](3,-1)(9,5)(9,5.4)(2.6,-1)
\pspolygon[linecolor=gray,fillcolor=gray,fillstyle=solid](-1,0)(9,0)(9,0.3)(-1,0.3)
\pspolygon[linecolor=gray,fillcolor=gray,fillstyle=solid](-1,2.5)(6,-1)(6.6,-1)(-1,2.8)
\psline[linewidth=1.5pt](-1,2.5)(6,-1)
\psline[linewidth=1.5pt](-1,0)(9,0)
\psline[linewidth=1.5pt](3,-1)(9,5)
\psline[linewidth=1.5pt](9,3.5)(2,7)
\psline[linewidth=1.5pt](-1,6)(9,6)
\psline[linewidth=1.5pt](-1,3.5)(6,7)
\psline[linewidth=1.5pt](0,-1)(0,7)
\rput(3.5,3.3){\small $\Delta^{(1)}$}
\end{pspicture}
\qquad \quad
\begin{pspicture}(-1,-1)(9,7)
\pspolygon(3,1)(4,1)(7,4)(5,5)(3,5)(1,4)(1,2)
\pspolygon[linewidth=1.5pt](4,0)(8,4)(4,6)(0,4)(0,2)
\rput(3.5,3.3){\small $\Delta^{(1)}$}
\rput(8,2){\small $\Delta^{\text{max}}$}
\end{pspicture}
\end{center}
When applying this construction to an \emph{arbitrary} two-dimensional lattice polygon $\Gamma$, one
ends up with a polygon $\Gamma^{(-1)}$ that is a lattice polygon if and only if $\Gamma = \Delta^{(1)}$ for some lattice polygon $\Delta$;
see \cite[Lem.\,10]{HaaseSchicho} for a proof of
this convenient criterion.
(If we call
\[ \left\{ \, \left. (i,j) \in \mathbb{R}^2 \, \right| \, <\! (i,j), v_\ell \!> \ = - a_\ell - 1 \, \right\} \]
the \emph{outward shift} of the edge corresponding to index $\ell$, then
a necessary, but generally insufficient condition
for $\Gamma^{(-1)}$ to be a lattice polygon is that the outward shifts of 
any pair of adjacent edges intersect
	in a lattice point \cite[Lem.\,9]{HaaseSchicho}.)\\
	
\noindent \emph{Remark.}
The criterion yields a method for algorithmically enumerating lattice polygons, as elaborated in \cite{movingout} and \cite[\S4.4]{Koelman}.
We will use this in the proofs of Theorem~\ref{maintheorem_gonality} and Theorem~\ref{maintheorem_neargonal}.\\

\noindent In Lemma~\ref{deltamaxcriterium} we will give a geometric interpretation of $\Delta^\text{max}$.

We use the notation $\mathcal{Z}(\cdot)$ to denote the algebraic set associated to an ideal, and $\mathcal{I}(\cdot)$
to denote the ideal of an algebraic set.

A \emph{curve} is always assumed irreducible, but we don't a priori
require it to be complete and/or smooth.
By the \emph{genus} of a curve $C$, which we denote by $g(C)$, we mean its geometric 
genus unless otherwise stated. The gonality of $C$ will be denoted by $\gamma(C)$.
A \emph{canonical curve} is a curve
that arises as the canonical image of a non-hyperelliptic smooth projective curve of genus $g \geq 3$. A \emph{canonical model}
of a curve $C$ is a canonical curve that is birationally equivalent to $C$.

\section{Divisors on toric surfaces} \label{section_toric}

This section gathers some facts on divisors on toric surfaces. 
Our primary objective is to fix notation and terminology, but we also group some statements
that are somewhat sprawled across our main references~\cite{coxlittleschenck,fulton}.

To a two-dimensional lattice polygon $\Delta$ we can associate a projective toric surface
$\text{Tor}(\Delta)$ over $k$, in two ways:
\begin{itemize}

 \item One can consider the (inner) normal fan $\Sigma_\Delta$, and let $\text{Tor}(\Delta) = \text{Tor}(\Sigma_\Delta)$ be
 the toric surface associated to it.

 \item One can define $\text{Tor}(\Delta)$ as the Zariski closure of the image of
 \begin{equation} \label{toricphi}
  \varphi_\Delta : \mathbb{T}^2 \hookrightarrow \mathbb{P}^N : (x,y) \mapsto
    \left( x^i y^j \right)_{(i,j) \in \Delta \cap \mathbb{Z}^2}
 \end{equation}
 (where $N =\sharp (\Delta \cap \mathbb{Z}^2) - 1$).
 Explicit equations for $\text{Tor}(\Delta)$ can be read from the combinatorics of $\Delta$, as follows.
 To each $(i,j) \in \Delta \cap \mathbb{Z}^2$ one associates a variable $X_{i,j}$.
 Then the ideal of $\text{Tor}(\Delta)$ is generated by the binomials
 \[ \prod_{\ell=1}^n X_{i_\ell,j_\ell} - \prod_{\ell=1}^n X_{i'_\ell,j'_\ell} \qquad \text{for which} \qquad \sum_{\ell=1}^n (i_\ell,j_\ell) = \sum_{\ell=1}^n (i'_\ell,j'_\ell) \]
 (apply \cite[Prop.\,2.1.4.(b,d)]{coxlittleschenck} to $\Delta \times \{1\} \subset \mathbb{R}^3$). A result of Koelman states that
 one can restrict to $n \in \{2,3\}$, and to $n=2$ as soon as $\partial \Delta \cap \mathbb{Z}^2 \geq 4$, see \cite{Koelman2,Schenck}.
\vspace{0.3cm}\\
 \noindent \emph{Examples.} \begin{itemize}
 \item $\text{Tor}(\Upsilon) = \mathcal{Z}(X_{0,0}^3 - X_{-1,-1}X_{1,0}X_{0,1}) \subset \mathbb{P}^3$,
 \item $\text{Tor}(\square) = \mathcal{Z}(X_{0,0}X_{1,1} - X_{1,0}X_{0,1}) \subset \mathbb{P}^3$,
 \item $\text{Tor}(\Gamma^5_1) = \mathcal{Z}(X_{0,0}^2 - X_{-1,0}X_{1,0}, X_{0,0}^2 - X_{0,-1}X_{0,1}) \subset \mathbb{P}^4$.
 \end{itemize}
\vspace{0.3cm}
\end{itemize}
Both constructions give rise to the same geometric object by \cite[Cor.\,2.2.19.(b)]{coxlittleschenck} and
the series of equivalences in the proof of \cite[Prop.\,6.1.10]{coxlittleschenck}. But
the second construction comes along with an embedding $\psi : \text{Tor}(\Delta) \hookrightarrow \mathbb{P}^N$, i.e.\
a very ample invertible sheaf $\psi^\ast \mathcal{O}_{\mathbb{P}^N}(1)$ on $\text{Tor}(\Delta)$.
Note that every complete fan in $\mathbb{R}^2$ arises as some $\Sigma_\Delta$. 

The self-action of $\mathbb{T}^2$
yields an action of $\mathbb{T}^2$ on $\varphi_\Delta(\mathbb{T}^2)$ that naturally extends to an action on all of $\text{Tor}(\Delta)$. The orbits of the latter
are in a dimension-preserving one-to-one correspondence with the faces of $\Delta$. Denote the Zariski closures
of the one-dimensional orbits (corresponding to the edges of $\Delta$ and to the rays of $\Sigma_\Delta$) by
$D_1, \dots, D_r$.
A Weil divisor that arises as a $\mathbb{Z}$-linear combination of the $D_\ell$'s is called \emph{torus-invariant}.
An important example is $K = - \sum_\ell D_\ell$, which is a canonical divisor; see \cite[Thm.\,8.2.3]{coxlittleschenck} or \cite[\S4.4]{fulton}.
To a torus-invariant Weil divisor $D = \sum_\ell a_\ell D_\ell$ one can associate the polygon
\begin{equation} \label{divisortopolygon}
\Delta_D = \bigcap_{\ell=1}^r H_\ell, \qquad \text{with $H_\ell = \left\{ \, \left. (i,j) \in \mathbb{R}^2 \, \right| \, <\! (i,j), v_\ell \! > \ \geq - a_\ell \, \right\}$,}
\end{equation}
where 
$v_\ell$ is the primitive generator of the corresponding ray in $\Sigma_\Delta$.
It can be proven \cite[Prop.\,4.3.3]{coxlittleschenck} that
\[ H^0(\text{Tor}(\Delta), D) = \left\{ \, \left. f \in k(x,y)^\ast \, \right| \, \text{div}(f) + D \geq 0 \, \right\} \cup \{0\} = \left\langle x^iy^j \right\rangle_{(i,j) \in \Delta_D \cap \mathbb{Z}^2} \]
(here $\langle \, \cdot \, \rangle$ denotes the $k$-linear span;
we view $x$ and $y$ as functions on $\text{Tor}(\Delta)$ through $\varphi_\Delta$).\\

\noindent \emph{Example.} Let $\Sigma$ be the fan given on the left, 
where the rays are enumerated as indicated.
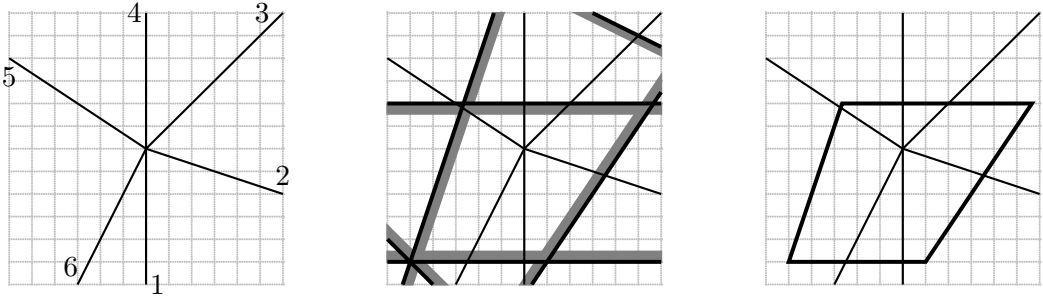
\begin{center}
\psset{unit=0.3cm}
\begin{pspicture}(-6,-6)(6,6.8)
\psgrid[subgriddiv=1,griddots=10,gridcolor=lightgray,gridlabels=0pt](-6,-6)(6,6)
\psline(0,-6)(0,6)
\psline(0,0)(6,-2)
\psline(0,0)(6,6)
\psline(0,0)(-6,4)
\psline(0,0)(-3,-6)
\rput(0.5,-6){\small $1$}
\rput(6,-1.2){\small $2$}
\rput(5.1,6){\small $3$}
\rput(-0.5,6){\small $4$}
\rput(-6,3.2){\small $5$}
\rput(-3.3,-5.2){\small $6$}
\end{pspicture}
\qquad \quad
\begin{pspicture}(-6,-6)(6,6.8)
\psgrid[subgriddiv=1,griddots=10,gridcolor=lightgray,gridlabels=0pt](-6,-6)(6,6)
\pspolygon[linecolor=gray,fillcolor=gray,fillstyle=solid](6,4.5)(3,6)(2.25,6)(6,4.125)
\pspolygon[linecolor=gray,fillcolor=gray,fillstyle=solid](6,-5)(-6,-5)(-6,-4.55)(6,-4.55)
\pspolygon[linecolor=gray,fillcolor=gray,fillstyle=solid](6,2)(-6,2)(-6,1.55)(6,1.55)
\pspolygon[linecolor=gray,fillcolor=gray,fillstyle=solid](-5.33,-6)(-1.33,6)(-0.9,6)(-4.9,-6)
\pspolygon[linecolor=gray,fillcolor=gray,fillstyle=solid](-4,-6)(-6,-4)(-6,-3.4)(-3.4,-6)
\pspolygon[linecolor=gray,fillcolor=gray,fillstyle=solid](6,3.25)(-0.17,-6)(0.33,-6)(6,2.5)
\psline[linewidth=1.5pt](6,4.5)(3,6)
\psline[linewidth=1.5pt](6,-5)(-6,-5)
\psline[linewidth=1.5pt](6,2)(-6,2)
\psline[linewidth=1.5pt](-4,-6)(-6,-4)
\psline[linewidth=1.5pt](0.33,-6)(6,2.5)
\psline[linewidth=1.5pt](-5.33,-6)(-1.33,6)
\psline(0,-6)(0,6)
\psline(0,0)(6,-2)
\psline(0,0)(6,6)
\psline(0,0)(-6,4)
\psline(0,0)(-3,-6)
\end{pspicture}
\qquad \quad
\begin{pspicture}(-6,-6)(6,6.8)
\psgrid[subgriddiv=1,griddots=10,gridcolor=lightgray,gridlabels=0pt](-6,-6)(6,6)
\pspolygon[linewidth=1.5pt](-5,-5)(1,-5)(5.66,2)(-2.66,2)
\psline(0,-6)(0,6)
\psline(0,0)(6,-2)
\psline(0,0)(6,6)
\psline(0,0)(-6,4)
\psline(0,0)(-3,-6)
\end{pspicture}
\end{center}
Let $D = 2D_1 + D_2 + 5D_3 + 5D_4 + D_5 + 3D_6$.
Then the corresponding half-planes are drawn in the middle, 
and
$\Delta_D$ is depicted on the right.
Remark that $\Delta_D$ is not a lattice polygon.
\\

One can also show that $D$ is Cartier if and only
if the apex of $H_\ell \cap H_m$ is an element of $\mathbb{Z}^2$ for each pair $\ell,m$ corresponding to adjacent 
edges of $\Delta$ \cite[Thm.\,4.2.8.(a,c)]{coxlittleschenck}.
If moreover every such apex is a vertex of $\Delta_D$ then $D$ is called \emph{convex} (in particular, if $D$ is a convex torus-invariant
Cartier divisor then $\Delta_D$ is a lattice polygon). If this gives a bijective apex-vertex correspondence
then $D$ is called \emph{strictly convex}.

A torus-invariant Cartier divisor $D$ is convex iff it is nef iff it is base-point free (i.e.\ $\mathcal{O}_{\text{Tor}(\Delta)}(D)$
is generated by global sections) by \cite[Thm.\,6.1.7~and~Thm.\,6.3.12]{coxlittleschenck}.
It is strictly convex iff it is ample iff it is very ample
\cite[Thm.\,6.1.14]{coxlittleschenck}. If $D$ is convex then all
higher cohomology spaces are trivial \cite[Thm.\,9.2.3]{coxlittleschenck}.
If $D_1$ and $D_2$ are convex torus-invariant Cartier divisors, 
then their intersection number
can be interpreted in terms of a mixed volume:
\[ D_1 \cdot D_2 = \text{MV}(\Delta_{D_1}, \Delta_{D_2}) = \text{Vol}(\Delta_{D_1} + \Delta_{D_2}) - \text{Vol}(\Delta_{D_1}) - \text{Vol}(\Delta_{D_2}), \]
where $\text{Vol}( \cdot )$ denotes the Euclidean area, and the addition of polygons is 
in Minkowski's sense (see \cite[\S5.3, first Cor.]{fulton} and the reasoning preceding \cite[\S5.5, (2)]{fulton}).
This is an instance of the Bernstein--Khovanskii--Koushnirenko (BKK) theorem.

Every Weil divisor
on $\text{Tor}(\Delta)$ is linearly equivalent to a
torus-invariant Weil divisor and two
equivalent torus-invariant Weil divisors $D_1$ and $D_2$
differ by some $\text{div}(x^iy^j)$  \cite[Thm.\,4.1.3]{coxlittleschenck}, so that the corresponding
polygons $\Delta_{D_1}$ and $\Delta_{D_2}$ are translates of each other.
Therefore, if one is willing to work modulo $\cong_t$, one can associate a polygon
$\Delta_D$ to \emph{any} Weil divisor $D$ (and a polygon $\Delta_\mathcal{L}$ to any invertible sheaf $\mathcal{L}$). All
definitions and statements above carry through.\\

\noindent \emph{Example.} 
We have $\Delta_{\psi^\ast \mathcal{O}_{\mathbb{P}^N}(1)} \cong_t \Delta$. Indeed, 
using (\ref{divisortopolygon}) it is straightforward to construct
a convex torus-invariant Cartier divisor $D_\Delta$ such that $\Delta_{D_\Delta} = \Delta$.
But then the global sections of $\mathcal{O}_{\text{Tor}(\Delta)}(D_\Delta)$
and $\psi^\ast \mathcal{O}_{\mathbb{P}^N}(1)$ are naturally identified. Since both sheaves
are globally generated, we find that
$\mathcal{O}_{\text{Tor}(\Delta)}(D_\Delta) \cong \psi^\ast \mathcal{O}_{\mathbb{P}^N}(1)$,
from which the claim follows.\\

\noindent \emph{Example.} Let $f \in k[x^{\pm 1}, y^{\pm 1}]$ be an irreducible Laurent polynomial and let $U(f)$ be the 
curve in $\mathbb{T}^2$ that it cuts out. Let $\Delta$ be any two-dimensional lattice polygon and let $C$ be the Zariski closure
of $\varphi_\Delta(U(f))$ in $\text{Tor}(\Delta)$. Let $P(f,\Delta)$ be the smallest convex polygon such that 
\begin{itemize}
  \item $\Delta(f) \subset P(f,\Delta)$, and 
  \item all edges of $P(f,\Delta)$
  are parallel to an edge of $\Delta$. 
\end{itemize}  
We claim that $\Delta_C \cong_t P(f,\Delta)$. Indeed, consider the torus-invariant Weil divisor $D_C = C - \text{div}(f)$,
so that we can assume that $\Delta_C = \Delta_{D_C}$. Then 
$f \in H^0(\text{Tor}(\Delta), D_C)$ and therefore $\Delta_C$ must contain
the support of $f$. Moreover, as we are working on $\text{Tor}(\Delta)$, every edge of $\Delta_C$ must be parallel to an edge of $\Delta$.
Each such edge must meet at least one point of the support of $f$, because otherwise the pole order of $f$ at the corresponding
torus-invariant prime divisor would be too large \cite[Prop.\,4.1.1]{coxlittleschenck}. So $\Delta_C$ must be the tightest fit, which is precisely $P(f,\Delta)$.


\section{Non-degenerate curves as smooth curves on toric surfaces} \label{section_nondegistoric}

We show how non-degenerate Laurent polynomials
naturally give rise to smooth curves in toric surfaces, and discuss how the non-degeneracy condition can be relaxed slightly. 
Much of the material below can be found
(possibly in disguised terms) in \cite{batyrev,CDV,coxlittleschenck}. On the other hand, Lemmata~\ref{deltamaxcriterium}--\ref{counterexample} seem genuinely new.

\subsection*{Non-degenerate curves.} Let $\Delta$ be a two-dimensional lattice polygon and consider a $\Delta$-non-degenerate Laurent polynomial
\[ f = \sum_{(i,j) \in \Delta \cap \mathbb{Z}^2} c_{i,j} x^iy^j \in k[x^{\pm 1}, y^{\pm 1}]. \]
Let $C \subset \text{Tor}(\Delta)$ be the Zariski closure of $\varphi_\Delta(U(f))$.
From the non-degeneracy of $f$ one sees that $C$ cuts out a smooth 
codimension $1$ subscheme in every $\mathbb{T}^2$-orbit of $\text{Tor}(\Delta)$. Because $\text{Tor}(\Delta)$ is 
normal \cite[Thm.\,3.1.5]{coxlittleschenck}, 
this is equivalent
to saying that $C$ is a smooth curve not containing any of
the zero-dimensional toric orbits and intersecting the
one-dimensional orbits transversally.
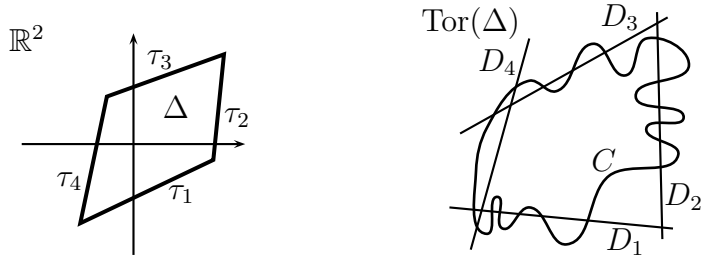
\begin{center}
\psset{unit=0.7cm}
\begin{pspicture}(-2.2,-2.2)(2.2,2.6)
\pspolygon[linewidth=1.5pt,linecolor=black](-1,-1.5)(1.5,-0.3)(1.7,1.7)(-0.5,0.9)
\psline{->}(-2.1,0)(2.1,0) \psline{->}(0,-2.1)(0,2.1)
\rput(0.8,0.7){$\Delta$} \rput(0.9,-0.95){$\tau_1$}
\rput(1.95,0.5){$\tau_2$} \rput(0.5,1.6){$\tau_3$}
\rput(-1.2,-0.7){$\tau_4$} \rput(-2,2){$\mathbb{R}^2$}
\end{pspicture} \qquad \qquad \qquad
\begin{pspicture}(-2.5,-2.2)(2.1,2.6)
\psccurve[linewidth=1pt](-1.8,-1.7)(-1.6,-1)(-1.5,-1.6)(-1,-1.2)(-0.3,-1.9)(0.5,-0.6)(1.8,-0.2)%
                        (1.2,0.3)(1.9,0.5)(1,0.9)(2,1.5)(1.2,2)(0.8,1.3)(0.2,1.9)(-0.5,1)(-1,1.2)%
                        (-1.8,0.3)(-2,-0.5)
\psline{-}(-2.5,-1.2)(1.7,-1.6) \psline{-}(1.5,-1.8)(1.4,2.5)
\psline{-}(1.6,2.4)(-2.3,0.2) \psline{-}(-1,2)(-2.1,-2)
\rput(0.4,-0.3){$C$} \rput(1.9,-1){$D_2$}
\rput(-1.6,1.55){$D_4$}
\rput(0.8,-1.85){$D_1$}
\rput(0.7,2.4){$D_3$}
\rput(-2.1,2.3){$\text{Tor}(\Delta)$}
\end{pspicture}
\end{center}
Note that $C$ is just the hyperplane section
\[ \sum_{(i,j) \in \Delta \cap \mathbb{Z}^2} c_{i,j}X_{i,j} \, = \, 0. \]
Therefore $\mathcal{O}_{\text{Tor}(\Delta)}(C) \cong \psi^\ast \mathcal{O}_{\mathbb{P}^N}(1)$. In particular, $C$
is a strictly convex Cartier divisor and $\Delta_C \cong_t \Delta$.

Toric surfaces are Cohen-Macaulay~\cite[Thm.\,9.2.9]{coxlittleschenck} and therefore
enjoy a nice adjunction theory, which we will use in the following form. 
Let $D_C$ be a torus-invariant
divisor that is linearly equivalent to $C$; for instance one may take $D_C = C - \text{div}(f)$. Then 
there is a canonical divisor $K_C$ on $C$ along
with an exact sequence
\begin{equation} \label{adjunctionsequence}
 0 \rightarrow \mathcal{O}_{\text{Tor}(\Delta)}(K) \rightarrow \mathcal{O}_{\text{Tor}(\Delta)}(D_C+K) \rightarrow \mathcal{O}_C(K_C) \rightarrow 0
\end{equation}
of morphisms of sheaves of $\mathcal{O}_{\text{Tor}(\Delta)}$-modules; locally the maps are given by $\cdot f$ and restriction to $C$, respectively. 

The existence of such an exact sequence is (in far greater generality) well-known to specialists in birational geometry; for example, 
this is essentially covered by
\cite[Prop.\,5.73]{kollar}. 
However we could not find a ready-to-use statement in the literature, so let us include the following flexible argument, which
was explained to us by Karl Schwede. 
Consider the short exact sequence
\[ 0 \rightarrow \mathcal{O}_{\text{Tor}(\Delta)}(-D_C) \stackrel{\cdot f}{\rightarrow} \mathcal{O}_{\text{Tor}(\Delta)} \rightarrow \mathcal{O}_{C} \rightarrow 0 \]
\vspace{-0.3cm}

\noindent and note that $\mathcal{O}_{\text{Tor}(\Delta)}(K)$ is a so-called dualizing sheaf for $\text{Tor}(\Delta)$; see \cite[\S4.4]{fulton}.
We apply the sheafy $\mathcal{H}om(\cdot, \mathcal{O}_{\text{Tor}(\Delta)}(K))$-functor to form a long exact sequence
\vspace{0.15cm}
\small \[  0 \rightarrow \mathcal{H}om(\mathcal{O}_C, \mathcal{O}_{\text{Tor}(\Delta)}(K)) \rightarrow \mathcal{H}om(\mathcal{O}_{\text{Tor}(\Delta)}, \mathcal{O}_{\text{Tor}(\Delta)}(K)) \rightarrow \mathcal{H}om(\mathcal{O}_{\text{Tor}(\Delta)}(-D_C), \mathcal{O}_{\text{Tor}(\Delta)}(K)) \qquad   \]
\[ \ \  \quad \qquad \qquad \qquad \qquad \qquad \qquad \qquad \qquad \rightarrow \mathcal{E}xt^1(\mathcal{O}_C, \mathcal{O}_{\text{Tor}(\Delta)}(K)) \rightarrow \mathcal{E}xt^1(\mathcal{O}_{\text{Tor}(\Delta)}, \mathcal{O}_{\text{Tor}(\Delta)}(K)).  \]
\normalsize
\vspace{-0.3cm}

\noindent The first term vanishes because $\mathcal{O}_C$ is torsion while $\mathcal{O}_{\text{Tor}(\Delta)}(K)$ is not.
The last term vanishes by~\cite[III.Prop.\,6.3(b)]{hartshorne}.
Finally because $\text{Tor}(\Delta)$ is Cohen-Macaulay, by \cite[Thm.\,2.12(1)]{reid} the fourth term is a dualizing sheaf
for $C$.
This is just $\mathcal{O}_C(K_C)$ and \eqref{adjunctionsequence} follows.

 Now note that $\Delta_K = \emptyset$, so that $H^0(\text{Tor}(\Delta),K) = 0$.  
Also $H^1(\text{Tor}(\Delta),K) = 0$, because by toric Serre duality \cite[Thm.\,9.2.10]{coxlittleschenck} the left-hand side is isomorphic to $H^1(\text{Tor}(\Delta),0)$, which vanishes by 
Demazure's theorem~\cite[Thm.\,9.2.3]{coxlittleschenck}.
Thus by taking the cohomology of (\ref{adjunctionsequence}) one finds that the restriction map
\begin{equation} \label{adjunctionisom}
 H^0(\text{Tor}(\Delta),D_C+K) \rightarrow H^0(C,K_C)
\end{equation}
is an isomorphism. Since the polygon associated to $D_C + K$ equals $\Delta^{(1)}$, we
recover
the well-known fact that
$g(C) = \sharp (\Delta^{(1)} \cap \mathbb{Z}^2)$. In fact, the isomorphism (\ref{adjunctionisom}) also
shows that 
\begin{equation} \label{canonicalembedding}
 \varphi_{\Delta^{(1)}}|_{U(f)} = 
\kappa \circ \varphi_\Delta|_{U(f)},
\end{equation}
where $\kappa : C \rightarrow \mathbb{P}^{g(C)-1}$ is a canonical morphism.
This seems less readily known, and will play an important role in Section~\ref{section_maroni}.
Using that the canonical image is rational iff $C$ is hyperelliptic, 
this observation implies the aforementioned fact
that $C$ is hyperelliptic iff the interior lattice points of $\Delta$ are collinear; see \cite[Lem.\,3.2.9]{Koelman} or \cite[Lem.\,2]{CaCo}
for more details.
If $C$ is non-hyperelliptic (i.e.\ $\Delta^{(1)}$ is two-dimensional) it follows
that the canonical image $\kappa(C)$ lies in $\text{Tor}(\Delta^{(1)}) \subset \mathbb{P}^{g(C)-1}$.\\

\noindent \emph{Remark.} If $C$ is an arbitrary (possibly singular, possibly non-Cartier) complete
curve on $\text{Tor}(\Delta)$ then the above adjunction process remains valid: one can still pick a torus-invariant
divisor $D_C$ that is equivalent to $C$, say with polygon $\Delta_C$ (not necessarily a lattice polygon!), and one will still find that the restriction map (\ref{adjunctionisom}) is an
isomorphism. When interpreting the outcome, some prudence is needed:
\begin{itemize}
\item In the non-Cartier case, note that in general
$\Delta_C^{(1)}$ is \emph{not} the polygon associated to $D_C+K$, which is the polygon obtained from $\Delta_C$ by shifting the edges inwards: this could result in
a polygon having vertices outside the lattice. But the lattice points of both polygons are the same, so in the smooth case 
it remains justified to say that $g(C) = \sharp (\Delta_C^{(1)} \cap \mathbb{Z}^2)$.
\item In the singular case 
%
we find that $\sharp (\Delta_C^{(1)} \cap \mathbb{Z}^2)$ is the \emph{arithmetic} genus of $C$, rather than
its geometric genus.\\
\end{itemize}

We note that classical adjunction theory, as elaborated in most textbooks, requires the ambient surface to
be smooth. Even though $\text{Tor}(\Delta)$ need not be smooth, it is possible to prove the genus formula
$g(C) = \sharp(\Delta^{(1)} \cap \mathbb{Z}^2)$ in this way, by first resolving the singularities using a toric blow-up. This is the
approach that is followed in~\cite[\S10.5]{coxlittleschenck}, for instance. We will briefly sketch this blow-up process
and show that it does not affect the combinatorics of $\Delta$,
because for the application
of Serrano's Theorem~\ref{thm_serrano} later on, we \emph{will} need that the ambient toric surface is smooth. (Serrano's
theorem plays the key role in the proofs of Theorems~\ref{maintheorem_gonality} and~\ref{maintheorem_neargonal}.)
So pick
a subdivision $\Sigma'$ of $\Sigma_\Delta$ such that the induced
birational morphism $\mu : \text{Tor}(\Sigma') \rightarrow \text{Tor}(\Sigma_\Delta)$ is a resolution
of singularities \cite[Thm.\,10.1.10]{coxlittleschenck}.
Let $C'$ be the strict transform of $C$ under $\mu$. By non-degeneracy $C'$ does not
meet the exceptional locus of $\mu$, so $C' = \mu^\ast C \cong C$. 
Note that $C'$ is again Cartier and convex, 
although not strictly convex (unless the subdivision is trivial).
It moreover remains true that $\Delta_{C'} \cong_t \Delta$. To prove this, we can suppose 
that $\Sigma'$ is obtained
from $\Sigma_\Delta$ by inserting a single ray $\sigma'$ (the general case then follows
by repeating the argument). Let $D_1, \dots, D_{r-2}$ be the torus-invariant prime divisors
on $\text{Tor}(\Sigma_\Delta)$ corresponding to the rays of $\Sigma_\Delta$ that are non-adjacent
to $\sigma'$, and let $D_1', D_2', \dots, D_{r-2}'$ be the according torus-invariant prime divisors
on $\text{Tor}(\Sigma')$. Then $D_i' = \mu^\ast D_i$ for all $i = 1, \dots, r-2$ (since $D_i'$ does not meet the exceptional 
locus of $\mu$).
Now by adding a divisor
of the form $\text{div}(x^iy^j)$ if needed, we see that
$C$ is linearly equivalent to a torus-invariant Weil divisor 
of the form $\sum_{\ell = 1}^{r-2} a_\ell D_\ell$.
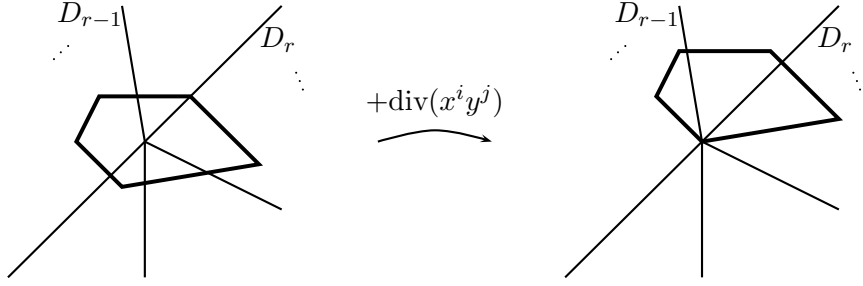
\begin{center}
\psset{unit=0.3cm}
\begin{pspicture}(-6,-6.4)(6,7.2)
\pspolygon[linewidth=1.5pt](-3,0)(-1,-2)(5,-1)(2,2)(-2,2)
\psline(0,0)(-1,6)
\psline(0,0)(-6,-6)
\psline(0,0)(0,-6)
\psline(0,0)(6,-3)
\psline(0,0)(6,6)
\rput(-2.4,5.6){\small $D_{r-1}$}
\rput(5.85,4.5){\small $D_r$}
\psarc[linestyle=dotted](0,0){7.3}{18}{25}
\psarc[linestyle=dotted](0,0){5.4}{128}{138}
\end{pspicture}
\qquad
\begin{pspicture}(-1,-6.4)(5,7.2)
\pscurve{->}(0,0)(2.5,0.5)(5,0)
\rput(2.5,1.8){\small $+ \text{div}(x^iy^j)$}
\end{pspicture}
\qquad
\begin{pspicture}(-6,-6.4)(6,7.2)
\pspolygon[linewidth=1.5pt](-2,2)(0,0)(6,1)(3,4)(-1,4)
\psline(0,0)(-1,6)
\psline(0,0)(-6,-6)
\psline(0,0)(0,-6)
\psline(0,0)(6,-3)
\psline(0,0)(6,6)
\rput(-2.4,5.6){\small $D_{r-1}$}
\rput(5.85,4.5){\small $D_r$}
\psarc[linestyle=dotted](0,0){7.3}{18}{25}
\psarc[linestyle=dotted](0,0){5.4}{128}{138}
\end{pspicture}
\end{center}
But then $C' = \mu^\ast C \sim \sum_{\ell = 1}^{r-2} a_\ell \mu^\ast D_\ell = \sum_{\ell = 1}^{r-2} a_\ell D'_\ell$,
from which it follows that $\Delta_{C'} \cong_t \Delta_C \cong_t \Delta$. 


\subsection*{$\Delta$-toric curves.} We now present a (slight) relaxation of
the non-degeneracy condition. Let $\Delta$ be a two-dimensional lattice polygon. We say that an irreducible Laurent
polynomial $f \in k[x^{\pm 1}, y^{\pm 1}]$ is \emph{$\Delta$-toric} if 
\begin{itemize}
  \item[(i)] $\Delta(f) \subset \Delta$,
  \item[(ii)] $\Delta(f)$ contains at least one point of every edge of $\Delta$, i.e.\ $P(f,\Delta) = \Delta$, and 
  \item[(iii)] the Zariski closure $C$ of $\varphi_\Delta(U(f))$ 
           is a smooth curve in $\text{Tor}(\Delta)$.
\end{itemize}
The condition that $P(f,\Delta) = \Delta$ ensures
that $C$ again arises as a hyperplane section of
$\text{Tor}(\Delta)$. We therefore still find that $\Delta_C \cong_t \Delta$.
All other conclusions of the preceding section remain valid, except for the part on resolutions
of singularities, where we add the assumption that
$\Sigma'$ does not subdivide any of the smooth cones of $\Sigma_\Delta$. Indeed, if it would, then 
this could affect $\Delta_C$. But since in practice there is no need for subdividing smooth cones, this is not an issue.
We also still obtain that $g(C) = \sharp (\Delta^{(1)} \cap \mathbb{Z}^2 )$ and that there is 
a canonical map $\kappa : C \rightarrow \mathbb{P}^{g(C)-1}$ satisfying 
(\ref{canonicalembedding}). Remark that H.\,Baker's bound \cite{beelen} implies
$g(C) \leq \sharp \left( \Delta(f)^{(1)} \cap \mathbb{Z}^2 \right)$, which together with $\Delta(f) \subset \Delta$
yields $\Delta(f)^{(1)} = \Delta^{(1)}$, a fact which
can also be proved directly by making a local analysis at the zero-dimensional $\mathbb{T}^2$-orbits of
$\text{Tor}(\Delta)$.

Geometrically, the only difference with $\Delta$-non-degeneracy
is that we allow $C$ to
contain some of the non-singular zero-dimensional orbits, or
to be tangent to some of the one-dimensional orbits. It cannot
pass through any of the singular zero-dimensional orbits however: otherwise $C$ would be singular as well. 

A curve $C/k$ is called \emph{$\Delta$-toric} if it is birationally equivalent to
$U(f)$ for a $\Delta$-toric Laurent polynomial $f$. This notion captures
all smooth projective curves on toric surfaces, as we will prove in Lemma~\ref{smoothmeanstoric} below (while
this is not true for non-degenerate curves: see Lemma~\ref{counterexample}).\\

\noindent \emph{Remark.} In the definition of being $\Delta$-toric, condition (iii) can be replaced by requiring
that 
\begin{itemize}
\item[(iii')] $g(C) = \sharp (\Delta^{(1)} \cap \mathbb{Z}^2)$. 
\end{itemize}
Indeed, in this case $C$ is automatically smooth,
because by adjunction theory $\sharp (\Delta^{(1)} \cap \mathbb{Z}^2)$ equals the arithmetic genus, which in the case of
singular curves is always strictly less than 
the geometric genus~\cite[IV.Ex.\,1.8]{hartshorne}. Recall that (iii') also
implies $\Delta(f)^{(1)} = \Delta^{(1)}$ by Baker's bound, which in turn implies (ii) as soon as $\Delta^{(1)}\neq \emptyset$.\\

Here is a geometric interpretation for the polygon $\Delta^\text{max}$ from Section~\ref{section_notation}.

\begin{lemma} \label{deltamaxcriterium}
  Let $\Delta$ be a lattice polygon and assume that $\Delta^{(1)}$ is two-dimensional. 
  Let $f \in k[x^{\pm 1}, y^{\pm 1}]$ be $\Delta$-toric and let $C$ be the Zariski closure of $\varphi_\Delta(U(f))$ in $\emph{Tor}(\Delta)$.
  Let $\kappa$ be as in (\ref{canonicalembedding}), so
  that $\kappa(C)$ can be viewed as a curve in the toric surface $\emph{Tor}(\Delta^{(1)})$. Then $\Delta_{\kappa(C)} \cong_t \Delta^\emph{max}$.
\end{lemma}

\noindent \textsc{Proof.} 
We see from (\ref{canonicalembedding})
that $\kappa(C)$ is the Zariski closure of $\varphi_{\Delta^{(1)}}(U(f))$ in $\text{Tor}(\Delta^{(1)})$.
From the remark concluding Section~\ref{section_toric} it follows that $\Delta_{\kappa(C)}$ is equivalent
to $P(f,\Delta^{(1)})$, the tightest polygon containing $\Delta(f)$ all of whose edges are parallel to an edge of $\Delta^{(1)}$.
But this polygon is clearly $\Delta^\text{max} = \Delta^{(1)(-1)}$.
\hfill $\blacksquare$\\

We now show that all smooth curves on toric surfaces are $\Delta$-toric, for an appropriate instance of $\Delta$.

\begin{lemma} \label{smoothmeanstoric}
  Let $C$ be a non-rational smooth projective curve on a toric surface, and
  let 
  \[ \tilde{\Delta}_C = \conv (\Delta_C \cap \mathbb{Z}^2). \]
  Then $C$ is $\tilde{\Delta}_C$-toric. 
\end{lemma}

\noindent Note that if $\Delta_C$ is a lattice polygon (i.e.\ if $C$ is Cartier) then $\tilde{\Delta}_C = \Delta_C$.
The non-rationality condition is not really a restriction: all smooth rational curves are isomorphic to $\mathbb{P}^1$, hence $\Sigma$-toric.\\

\noindent \textsc{Proof.} 
Let $X$ be our toric surface, containing the torus $\mathbb{T}^2$ as an open subset. 
Since $C$ is non-rational, it is non-torus-invariant. So $C \cap \mathbb{T}^2$ is
defined by an irreducible Laurent polynomial $f \in k[x^{\pm 1}, y^{\pm 1}]$.
The torus-invariant divisor $D_C = C - \text{div}(f)$ is equivalent to $C$, so that we can assume that $\Delta_C$ is the polygon associated
to $D_C$. Because $f \in H^0(X, D_C)$ we see that $f$ is supported on $\Delta_C$, and because $\Delta(f)$ is a lattice polygon we even have that
\[  \Delta(f) \subset \tilde{\Delta}_C \subset \Delta_C \]
and in particular that
\begin{equation} \label{tildeinclusion}
  \sharp(\Delta(f)^{(1)} \cap \mathbb{Z}^2) \leq \sharp(\tilde{\Delta}^{(1)}_C \cap \mathbb{Z}^2) \leq \sharp (\Delta^{(1)}_C \cap \mathbb{Z}^2).
\end{equation}
By adjunction theory the genus of $C$ equals $\sharp (\Delta_C^{(1)} \cap \mathbb{Z}^2)$. On the other
hand by Baker's bound it is at most $\sharp (\Delta(f)^{(1)} \cap \mathbb{Z}^2)$. Thus the inequalities in \eqref{tildeinclusion}
are equalities, and in particular the genus of $C$ also equals $\sharp (\tilde{\Delta}_C^{(1)} \cap \mathbb{Z}^2)$. In other words, with respect to the lattice polygon
$\tilde{\Delta}_C$, our polynomial $f$ satisfies
condition (iii') mentioned above, and therefore it is $\tilde{\Delta}_C$-toric. \hfill $\blacksquare$\\

From the proof we see that $C$ is in fact also $\Delta(f)$-toric, but we chose to provide a polygon
that depends on the divisor class of $C$ only (up to translation).
As an immediate corollary to the previous lemmata and their proofs, we obtain:

\begin{lemma} \label{cruciallemma}
 Let $\Delta$ be a lattice polygon and assume that $\Delta^{(1)}$ is two-dimensional. 
  Let $f \in k[x^{\pm 1}, y^{\pm 1}]$ be $\Delta$-toric.
Then $f$ is also $\Delta^\emph{max}$-toric.
\end{lemma}

\noindent This lemma
will play an important role in the proofs
of Theorems~\ref{maintheorem_gonality} and~\ref{maintheorem_neargonal}. 
It is in the same vein as Kawaguchi's notion of 
\emph{relative minimality} \cite[Def.\,3.9]{Kawaguchi}, and can
be proven more directly,
by noting that $\Delta$ is obtained from $\Delta^\text{max}$ by clipping off
a number of vertices, without affecting the interior. From Pick's theorem
it follows that such a vertex is necessarily smooth, i.e.\ that the primitive normal vectors
of its adjacent edges form a basis of $\mathbb{Z}^2$. Then locally around the corresponding
zero-dimensional orbit, $\text{Tor}(\Delta^\text{max})$ looks like $\mathbb{A}^2$ with $C$ 
passing smoothly through the origin. The smoothness of $C$ outside these zero-dimensional orbits
then just follows from the fact that $f$ is $\Delta$-toric.
\\

Clearly every $\Delta$-non-degenerate curve is $\Delta$-toric. The converse implication
may fail:

\begin{lemma} \label{counterexample}
There exist instances of two-dimensional lattice polygons $\Delta$, along with $\Delta$-toric curves that
are not $\Delta$-non-degenerate. More precisely,
let 
\[ f = 1 + x^5 + y^2 + x^2y^3 \in k[x^{\pm 1}, y^{\pm 1}] \quad \text{and} \quad \Delta = \conv\{ (0,0), (5,0), (2,3), (0,3) \}. \]
Then $f$ is $\Delta$-toric, but $U(f)$ is not $\Delta$-non-degenerate, that is, it is not birationally equivalent to $U(f')$ for
some $\Delta$-non-degenerate Laurent polynomial $f' \in k[x^{\pm 1},y^{\pm 1}]$
\end{lemma}

\noindent \textsc{Proof.}
Our proof uses the theory of trigonal curves. We need the following facts. If $C / k$ is a trigonal curve of genus $g \geq 5$,
then the intersection of all quadrics
  containing its canonical model $C_\text{can} \subset \mathbb{P}^{g-1}$ is a rational normal surface scroll $S$ spanned by 
  two rational normal curves $R_1$ and $R_2$ of respective degrees $e_1$ and $e_2$, where $e_1 \leq e_2$. These numbers are uniquely
  determined and are called the Maroni
  invariants of $C_\text{can}$. See \cite[(4.11)]{saintdonat} for a proof, and \cite[Ex.\,8.17]{harris} and 
  Section~\ref{section_maroni} for more background on this terminology. 
For our needs it is important that if $e_1 < e_2$ then $R_1$ is uniquely determined by $S$ \cite[Prop.\,8.20(b)]{harris}.
It follows that in the case where $e_1 < e_2$, the number of points at which $C_\text{can}$ is tangent to $R_1$ is an invariant of $C$, which
we denote by $t_C$.


Now the reader can verify that $f$ is indeed $\Delta$-toric, i.e.\ the Zariski
closure $C$
of $\varphi_\Delta(U(f))$ is a smooth curve in $\text{Tor}(\Delta)$. Note that $C$ is a trigonal curve of genus $5$, since 
it is non-hyperelliptic by \cite[Lem.\,3.2.9]{Koelman} and
the map $U(f) \rightarrow \mathbb{T}^1 : 
(x,y) \mapsto x$ is of degree $3$.
Let $C_\text{can}$ be the canonical model obtained by taking the
Zariski closure of $\varphi_{\Delta^{(1)}}(U(f))$ inside $\text{Tor}(\Delta^{(1)}) \subset \mathbb{P}^4$. 
Since the latter surface 
is generated by quadrics, it must be our rational normal scroll $S$. 
The scrollar structure can easily be made explicit in this case. 
In particular, one verifies that $e_1 = 1$ and $e_2 = 2$, and that the line $R_1$ is the torus-invariant
prime divisor of $\text{Tor}(\Delta^{(1)})$ corresponding to the top edge of $\Delta^{(1)}$.
Now remark that $\Sigma_\Delta = \Sigma_{\Delta^{(1)}}$, so we have a natural isomorphism
$\mu: \text{Tor}(\Delta^{(1)}) \rightarrow \text{Tor}(\Delta)$, which is compatible with 
the respective embeddings of $\mathbb{T}^2$ in $\text{Tor}(\Delta^{(1)})$ and 
$\text{Tor}(\Delta)$, i.e.\ $\varphi_\Delta = \mu \circ \varphi_{\Delta^{(1)}}$. In particular
$\mu(C_\text{can}) = C$, and because $\mu$ behaves well with respect to the toric orbits we find that $t_C$ can
be interpreted as the number of points at which $C$ is tangent to the torus-invariant prime divisor of $\text{Tor}(\Delta)$ 
corresponding to the top edge of $\Delta$. Using this, one easily checks that
$t_C = 1$. On the other hand, the same reasoning shows that
if $U(f)$ were $\Delta$-non-degenerate, then $t_C$ would be $0$. \hfill $\blacksquare$\\

\noindent \emph{Remarks.} 
\begin{itemize}
\item It is not possible to construct similar counterexamples from arbitrary two-dimen\-sional lattice polygons.
For instance, let $\Delta = d\Sigma$ for some integer $d \geq 1$, so that $\text{Tor}(\Delta) \cong \mathbb{P}^2$. 
Then every $\Delta$-toric curve is $\Delta$-non-degenerate.
Indeed, using an automorphism of $\mathbb{P}^2$, every smooth projective plane curve can be positioned in such
a way that it does not contain any of the coordinate points, and such that it intersects the coordinate axes transversally.
\item In all theorems and lemmata appearing in Sections~\ref{section_gonality} to~\ref{section_maroni} of 
this paper (which contain our main results), the notions of being $\Delta$-non-degenerate and $\Delta$-toric are interchangeable, i.e.\ only the property of
being $\Delta$-toric is used in the proofs. For instance:
\begin{quote}
\noindent \textbf{Corollary~\ref{cacoconjecture}.}~\emph{Let 
$f \in k[x^{\pm 1},y^{\pm 1}]$ be non-degenerate with respect to its Newton polygon
 $\Delta = \Delta(f)$. Then the gonality of $U(f)$ equals $\emph{lw}(\Delta^{(1)}) + 2$,
 unless $\Delta^{(1)} \cong \Upsilon$ (i.e.\ $\Delta \cong 2\Upsilon$), in which case it equals $3$.}
\end{quote}
from Section~\ref{section_gonality} can be replaced by the slightly stronger statement 
that the gonality of a $\Delta$-toric curve equals $\text{lw}(\Delta^{(1)}) + 2$, unless
$\Delta^{(1)} \cong \Upsilon$, in which case it equals $3$. We have chosen to state our main results in a toric-geometry-free language, however.
\end{itemize}

\section{Lattice directions and combinatorial pencils} \label{section_combpencils}

\subsection*{Lattice directions}

A \emph{lattice direction} is just a primitive element of $\mathbb{Z}^2$.
For a non-empty lattice polygon $\Delta$ and a lattice direction $v = (a,b)$, the \emph{width of $\Delta$ with respect to $v$} is
the minimal $d$ for which there exists an $m \in \mathbb{Z}$ such that $\Delta$ is contained in the strip
\[ m \leq aY - bX \leq m + d.\]
Note that $w(\Delta,v) = w(\Delta,-v)$.
If $w(\Delta,v) = d$, we will sometimes say that \emph{$v$ computes $d$}. It is convenient to define $w(\emptyset,v) = -1$.
(This notion appeared in \cite[Def.\,5]{LubbesSchicho} where it is called the \emph{viewangle width}.)\\

\noindent \emph{Example.} The width
of $d\Sigma$ with respect to $(1,1)$ is $2d$, while its width with respect to $(1,-1)$ is $d$.
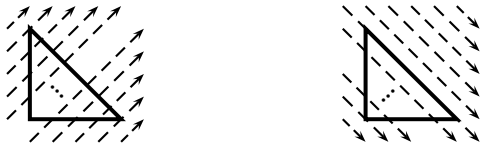
\begin{center}
\psset{unit=0.3cm}
\begin{pspicture}(-1,-1)(5,5)
\pspolygon[linewidth=1.5pt](0,0)(4,0)(0,4)
\psline[linestyle=dashed]{->}(-1,0)(4,5)
\psline[linestyle=dashed]{->}(-1,1)(3,5)
\psline[linestyle=dashed]{->}(-1,2)(2,5)
\psline[linestyle=dashed]{->}(-1,3)(1,5)
\psline[linestyle=dashed]{->}(-1,4)(0,5)
\psline[linestyle=dashed]{->}(0,-1)(5,4)
\psline[linestyle=dashed]{->}(1,-1)(5,3)
\psline[linestyle=dashed]{->}(2,-1)(5,2)
\psline[linestyle=dashed]{->}(3,-1)(5,1)
\psline[linestyle=dashed]{->}(4,-1)(5,0)
\rput(1.2,1.2){$\cdot$}
\rput(1,1.4){$\cdot$}
\rput(1.4,1){$\cdot$}
\end{pspicture}
\qquad \qquad \qquad
\begin{pspicture}(-1,-1)(5,5)
\pspolygon[linewidth=1.5pt](0,0)(4,0)(0,4)
\psline[linestyle=dashed]{->}(4,5)(5,4)
\psline[linestyle=dashed]{->}(3,5)(5,3)
\psline[linestyle=dashed]{->}(2,5)(5,2)
\psline[linestyle=dashed]{->}(1,5)(5,1)
\psline[linestyle=dashed]{->}(0,5)(5,0)
\psline[linestyle=dashed]{->}(-1,5)(5,-1)
\psline[linestyle=dashed]{->}(-1,4)(4,-1)
\psline[linestyle=dashed]{->}(-1,2)(2,-1)
\psline[linestyle=dashed]{->}(-1,1)(1,-1)
\psline[linestyle=dashed]{->}(-1,0)(0,-1)
\rput(1,1){$\cdot$}
\rput(0.8,0.8){$\cdot$}
\rput(1.2,1.2){$\cdot$}
\end{pspicture}
\end{center}

\begin{lemma} \label{numberofdirections}
If $\Delta$ is two-dimensional, then for a given $d \in \mathbb{Z}_{\geq 0}$, the number of lattice directions computing $d$ is finite.
\end{lemma}

\noindent \textsc{Proof.} It suffices to prove that for each $d$, the number of lattice directions $v$ for which
$w(\Delta,v) \leq d$
is finite.
Since $\Delta$ is two-dimensional we may assume that it contains
the standard simplex $\Sigma$ (see e.g.\ \cite[Prop.\,1.2.4.(b)]{Brunsetal}, although this easily follows from Pick's theorem), so that $w(\Sigma, v) \leq w(\Delta,v)$ for every $v$. Thus it suffices
to prove that for each $d$, the number of lattice directions $v$ for which $w(\Sigma,v) \leq d$ is finite.
But this is straightforward.
\hfill $\blacksquare$\\

Assume that $w(\Delta,v) = d \geq 2$. Write $v = (a,b)$ and 
assume that $\Delta$ is contained in the strip 
$m \leq aY - bX \leq m + d$.
Then we define the \emph{width invariants of $\Delta$ with respect to $v$} as the tuple
\[ E(\Delta,v) \ = \ \left( E_\ell \right)_{\ell = 1, \dots, d - 1} \]
where
\[ E_\ell = \sharp \left\{ \left. \, (i,j) \in \Delta^{(1)} \cap \mathbb{Z}^2 \, \right| \,  aj - bi = m + \ell \, \right\} - 1.\]
(The reason for the $-1$ term will become clear in Section~\ref{section_maroni}.)\\

\noindent \emph{Example.} Let $v = (1,0)$ and $d \in \mathbb{Z}_{\geq 2}$. Then $w(d\Sigma,v) = d$ and
$E(d\Sigma,v) = (d-3,d-4,d-5, \dots, 1,0,-1)$.\\

The \emph{lattice width of $\Delta$} is
\[ \text{lw}(\Delta) = \min_v w(\Delta,v).\]
Equivalently, $\text{lw}(\Delta)$ is the minimal $d$ such that $\Delta$ is unimodularly equivalent to a lattice polygon that is contained in
a horizontal strip of height $d$; for the latter, two lattice directions computing the lattice width are $(\pm 1, 0)$. 
If a lattice direction computes the lattice width, we call it a \emph{lattice width direction for $\Delta$}.\\

\noindent \emph{Example.} Let $d \in \mathbb{Z}_{\geq 0}$. Then $\text{lw}(d\Sigma) = d$. Indeed, clearly $\text{lw}(d\Sigma) \leq d$, while
$\text{lw}(\Delta) \geq d$ follows from the fact that every edge of $d\Sigma$ contains $d + 1$ lattice points.\\

A convenient tool for computing $\text{lw}(\Delta)$ is given by (i) from Lemma~\ref{latticewidthproperties}
below, which gathers some useful facts about the lattice width:
  
\begin{lemma} \label{latticewidthproperties}
Let $\Delta$ be a two-dimensional lattice polygon. 
\begin{enumerate}
  \item[(i)] One has $\emph{lw}(\Delta^{(1)}) = \emph{lw}(\Delta) - 2$, unless $\Delta \cong d\Sigma$ for some $d \geq 2$ in which case
        $\emph{lw}(\Delta^{(1)}) = \emph{lw}(\Delta) - 3 = d - 3$.
  \item[(ii)] A 
        lattice width direction for $\Delta$ is also a lattice width direction for $\Delta^{(1)}$;
        if moreover $\Delta^{(1)} \neq \emptyset$ and $\Delta^{(1)} \not \cong (d-3)\Sigma$ for any $d \geq 3$, then 
        the converse holds as well.
  \item[(iii)] Assume $\emph{lw}(\Delta) \geq 2$ and $\Delta \not \cong d\Sigma$ for any $d \geq 2$. Then the width invariants of $\Delta$ with respect to 
  a lattice width direction 
  are all non-negative.
  \item[(iv)] There are at most $4$ pairs $\pm v$ of lattice width directions for $\Delta$; the bound is met if and only
  if $\Delta \cong d \Gamma^5_1$ for some $d \in \mathbb{Z}_{\geq 1}$.
  \item[(v)] If $v_1,v_2$ are lattice width directions for $\Delta$, then $\left| \det(v_1,v_2) \right| \leq 2$;
  if equality holds then $\Delta \cong d \Gamma^5_1$ for some $d \in \mathbb{Z}_{\geq 1}$.
  \item[(vi)] One has $\emph{lw}(\Delta)^2 \leq \frac{8}{3} \emph{Vol}(\Delta)$, and equality holds if and only if $\Delta \cong d\Upsilon$ for some $d \geq 1$.
\end{enumerate}
\end{lemma}
\noindent \textsc{Proof.} For (i) and (ii) see \cite[Thm.\,4]{CaCo} or \cite[Thm.\,13]{LubbesSchicho}.

Claim (iii) can be proved by induction, as follows. Let $v = (a,b)$ be a lattice width direction for $\Delta$ and
let $m$ be such that $\Delta$ is contained in the strip $m \leq aY - bX \leq m + \text{lw}(\Delta)$.
We have to show that for each $\ell = 1, \dots, \text{lw}(\Delta) - 1$ there exists 
an $(i,j) \in \Delta^{(1)} \cap \mathbb{Z}^2$ 
such that $aj-bi = m + \ell$. Now (i) implies that this must be the case for
$\ell = 1$ and $\ell = \text{lw}(\Delta) - 1$, while from (ii) it follows that $v$ is also a lattice width direction for
$\Delta^{(1)}$. So the claim follows by recursively applying it
to $\Delta^{(1)}$; if at some point $\Delta^{(1)}$ happens to be of the form $d\Sigma$ for some $d \geq 2$ then
the claim can be verified explicitly.

For (iv) see \cite{draismaetal}. 

To prove (v), let $v_1,v_2$ be lattice width directions for $\Delta$ for which 
$\left| \det(v_1,v_2) \right| > 1$. Using a unimodular transformation if needed, we can assume
that $v_1 = (1,-1)$ and $v_2 = (a,b)$ with $a,b>0$, and that $\Delta$ is contained in the strips
$0 \leq Y + X \leq \text{lw}(\Delta)$ and $0 \leq aY - bX \leq \text{lw}(\Delta)$.
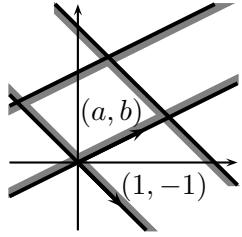
\begin{center}
\psset{unit=0.3cm}
\begin{pspicture}(-3,-3)(7,7)
\pspolygon[linecolor=gray,fillcolor=gray,fillstyle=solid](-3,3)(3,-3)(3.4,-3)(-3,3.4)
\pspolygon[linecolor=gray,fillcolor=gray,fillstyle=solid](-1,7)(7,-1)(6.6,-1)(-1,6.6)
\pspolygon[linecolor=gray,fillcolor=gray,fillstyle=solid](-3,-1.5)(7,3.5)(7,3.85)(-3,-1.15)
\pspolygon[linecolor=gray,fillcolor=gray,fillstyle=solid](-3,2.5)(6,7)(6.7,7)(-3,2.15)
\psline[linewidth=1.5pt](-3,3)(3,-3)
\psline[linewidth=1.5pt]{->}(0,0)(2,-2)
\psline[linewidth=1.5pt](-1,7)(7,-1)
\psline[linewidth=1.5pt](-3,-1.5)(7,3.5)
\psline[linewidth=1.5pt]{->}(0,0)(3,1.5)
\psline[linewidth=1.5pt](-3,2.5)(6,7)
\rput(3.8,-1.1){\small $(1,-1)$}
\rput(1.5,2.2){\small $(a,b)$}
\psline{->}(-3,0)(7,0)
\psline{->}(0,-3)(0,7)
\end{pspicture}
\end{center}
Thus $\Delta$ is contained in the parallelogram
\small
\[ \conv \left\{ (0,0), \left( \frac{a \text{lw}(\Delta)}{a+b}, \frac{b \text{lw}(\Delta)}{a+b} \right),
\left(\frac{(a-1)\text{lw}(\Delta)}{a+b}, \frac{(b+1)\text{lw}(\Delta)}{a+b} \right) , \left(- \frac{\text{lw}(\Delta)}{a+b}, \frac{\text{lw}(\Delta)}{a+b} \right) \right\}. \]
\normalsize
The horizontal width of this parallelogram equals $(a+1)\text{lw}(\Delta)/(a+b)$, while its vertical width equals
$(b+1)\text{lw}(\Delta)/(a+b)$. By the definition of $\text{lw}(\Delta)$ it follows
that $a=b=1$, so that $\left| \det(v_1,v_2) \right| = 2$. Moreover, these four vertices must be actual vertices of $\Delta$.
In particular, they must be contained in $\mathbb{Z}^2$, from which
one sees that $\text{lw}(\Delta)$ is even, and $\Delta \cong \frac{\text{lw}(\Delta)}{2}\Gamma^5_1$.

For (vi)\ see \cite{FejesMakai}. \hfill $\blacksquare$\\

Note that Lemma~\ref{latticewidthproperties}.(v) implies
that if $\Delta$ has two linearly independent lattice width directions $v_1, v_2$, then there is
a unimodular transformation mapping $\Delta$
inside $\text{lw}(\Delta) \square$. (Indeed, if $\left|\det (v_1, v_2) \right| = 1$ 
then one can take a $\mathbb{Z}$-linear transformation mapping $v_1$ to $(1,0)$ and $v_2$ to $(0,1)$, and compose
it with an appropriate translation; if $\left| \det (v_1, v_2) \right| \neq 1$ then $\Delta$ is of the form $d\Gamma^5_1$, 
and the statement can be verified explicitly.) 
In particular, it follows that
\begin{equation} \label{genussqbound}
\sharp (\Delta^{(1)} \cap \mathbb{Z}^2) \leq (\text{lw}(\Delta)- 1)^2
\end{equation}
in this case.\\

\noindent \emph{Example.} Let $\Delta$ be the lattice polygon
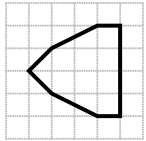
\begin{center}
\psset{unit=0.3cm}
\begin{pspicture}(-1,-3)(5,3)
\psgrid[subgriddiv=1,griddots=10,gridcolor=lightgray,gridlabels=0pt](-1,-3)(5,3)
\pspolygon[linewidth=1.5pt](0,0)(1,-1)(3,-2)(4,-2)(4,2)(3,2)(1,1)
\end{pspicture}
\end{center}
for which $\text{lw}(\Delta) = 4$ (as can be seen by applying Lemma~\ref{latticewidthproperties}.(i)).
Clearly $\pm (1,0)$ and $\pm (0,1)$ are lattice directions computing $\text{lw}(\Delta)$. It is also immediate
that $\Delta \not \cong d\Gamma^5_1$ for any $d \in \mathbb{Z}_{\geq 1}$, so that 
by Lemma~\ref{latticewidthproperties}.(iv)\ the number of pairs $\pm v$ of lattice width directions is either two or three.
But three is impossible, because by Lemma~\ref{latticewidthproperties}.(v)\ the
third pair would need to be among $\pm (1,1), \pm (1,-1)$, both of which correspond
to widths that strictly exceed $4$.\\

\noindent \emph{Remark.} Lemma~\ref{latticewidthproperties}.(iii) can also be proven 
using the well-known geometric fact that gonality pencils are always complete, by combining 
Theorem~\ref{maintheorem_gonality} and Corollary~\ref{nonnegativewidth} below.

\subsection*{Combinatorial pencils}

\noindent Returning to the geometric side, 
let $\Delta$ be a two-dimensional lattice polygon, let $v = (a,b)$ be a lattice direction,
let $f \in k[x^{\pm 1}, y^{\pm 1}]$ be a $\Delta$-non-degenerate or $\Delta$-toric Laurent polynomial, and
let $C \subset \text{Tor}(\Delta)$ be the corresponding
smooth projective curve (i.e.\ the Zariski closure of $\varphi_\Delta(U(f))$, where $\varphi_\Delta$
is as in Section~\ref{section_toric}).
We associate to $v$ a linear pencil $g_v$ on $C$ as follows. For each 
$c \in \mathbb{P}^1 = \mathbb{T}^1 \cup \{0, \infty \}$ we have a function
$x^ay^b - c$ on $\text{Tor}(\Delta)$ (where $x^ay^b - \infty$ should be read as $x^{-a}y^{-b}$) whose zero divisor
$\mathcal{F}_c$ cuts out a divisor $D_c$ on $C$. Then 
\[ g_v = \{D_c\}_{\substack{c \in \mathbb{P}^1 \\ \mathcal{F}_c \neq C}}.\]
In other words this is the trace on $C$ of the linear system $\{ \mathcal{F}_c \}_c$, in the sense of \cite[p.\,158]{hartshorne}.\\

\noindent \emph{Remark.} The subscript $\mathcal{F}_c \neq C$ is usually superfluous, but it could indeed happen that $\mathcal{F}_c = C$ for some $c$.
Example: $f = x + 1$, $\Delta = \Sigma$, $v = (1,0)$ and $c = -1$. In this example $g_v$ is just the linear system 
consisting of one base point, namely the point $(0:1:0) \in \mathbb{P}^2 = \text{Tor}(\Sigma)$ (note the abuse of language here).
By genus considerations $\mathcal{F}_c = C$ can occur only if $\Delta^{(1)} = \emptyset$. Since from Section~\ref{section_gonality} 
on, all our theorems and lemmata that involve combinatorial pencils exclude the case $\Delta^{(1)} = \emptyset$ (for other reasons), the reader 
can in fact ignore the possibility of this event.\\

There are several ways of seeing that $g_v$ has degree $w(\Delta,v)$.
One approach, the details of which we leave to the reader, uses the BKK theorem along with
the fact that $w(\Delta,v) = \text{MV}(\Delta, \conv \{ 0, v \})$. We will give a more elementary argument that
gives finer information.

\begin{lemma} \label{gv_basepointfree}
The pencil $g_v$ is of degree $w(\Delta,v)$. More precisely, it splits into a base-point free part of degree
$w(\Delta(f),v)$ and a fixed part of degree $w(\Delta,v) - w(\Delta(f),v)$ that is supported on at most two points.
In particular, if
\begin{itemize}
  \item $f$ is $\Delta$-non-degenerate, or
  \item $v$ is a lattice width direction, $\emph{lw}(\Delta) \geq 2$, and $\Delta \not \cong d\Sigma$ for all $d \geq 2$,
\end{itemize}
then $g_v$ is base-point free.
\end{lemma}

\noindent \textsc{Proof.} 
We will assume that $\mathcal{F}_c \neq C$ for all $c \in \mathbb{P}^1$, and leave the details of the other case to the reader.
Then the rational map $U(f) \rightarrow \mathbb{T}^1 : (x,y) \mapsto x^ay^b$ 
extends to a degree $w(\Delta(f),v)$ morphism $C \rightarrow \mathbb{P}^1$. Its fibers
determine a
base-point free pencil that necessarily matches with the base-point free part of $g_v$.

As for the fixed part, suppose that $\Delta$ and $\Delta(f)$ are contained in the strips 
\[ m \leq aY - bX \leq m + w(\Delta,v) \qquad \text{and} \qquad m_f \leq aY - bX \leq m_f + w(\Delta(f),v), \]
respectively.
If $\Delta$ has a unique vertex $v_\text{low}$ lying on the line $m = aY - bX$, the corresponding
zero-dimensional orbit $O(v_\text{low})$ is contained in every divisor $\mathcal{F}_c$. Similarly, if there
is a unique vertex $v_\text{top}$ on the line $m + w(\Delta,v) = aY - bX$ then $O(v_\text{top})$ is contained in every $\mathcal{F}_c$.
All other points of $\text{Tor}(\Delta)$ lie on a unique $\mathcal{F}_c$. This means that the fixed part of $g_v$
is supported on at most these two points. 
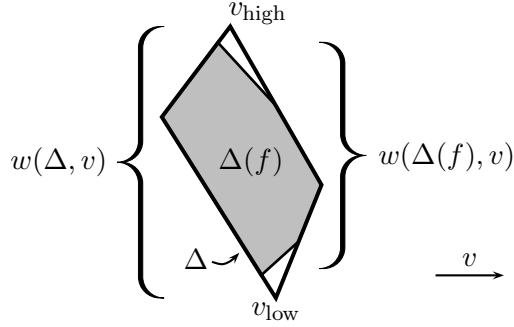
\begin{center}
\psset{unit=0.3cm}
\begin{pspicture}(-6,-7)(6,7.5)
\pspolygon[linewidth=1pt,fillstyle=solid,fillcolor=lightgray](4,-1)(2,2.5)(-0.55,5.3)(-3,2)(1.35,-5)(3,-3.5)
\pspolygon[linewidth=1.5pt](2,-6)(4,-1)(0,6)(-3,2)
\psline{->}(9,-5)(12,-5)
\rput(10.5,-4.3){\small $v$}
\psbrace(-2.9,6)(-2.9,-6){ }
\psbrace(3.9,-4.7)(3.9,5.3){ }
\rput(-7.5,0){\small $w(\Delta,v)$}
\rput(9.5,0.25){\small $w(\Delta(f),v)$}
\rput(1.2,6.6){\small $v_\text{high}$}
\rput(2,-6.6){\small $v_\text{low}$}
\rput(-1.5,-4.3){\small $\Delta$}
\pscurve{->}(-0.8,-4.4)(-0.3,-4.5)(0.4,-4)
\rput(0.9,0){\small $\Delta(f)$}
\end{pspicture}
\end{center}
Now if there is indeed a unique lower-most vertex $v_\text{low}$ of $\Delta$, then a local analysis shows
that a generic $\mathcal{F}_c$ intersects $C$ in $O(v_\text{low})$ with multiplicity $m_f - m$, or in other words, the order
of the fixed part of $g_v$ at $O(v_\text{low})$ equals $m_f - m$. 
If there is no unique lower-most vertex, then
necessarily $m = m_f$, otherwise there would be an edge of $\Delta$ not supporting any term of $f$, contradicting that $f$
is $\Delta$-toric. 
A similar analysis at the top then yields that the fixed part of $g_v$ has degree $w(\Delta,v) - w(\Delta(f),v)$.

For the last claim it suffices to note that if
$f$ is $\Delta$-non-degenerate then $\Delta(f) = \Delta$, and therefore $w(\Delta(f),v) = w(\Delta,f)$, while
if $v$ is a lattice width direction, $\text{lw}(\Delta) \geq 2$, and $\Delta \not \cong d\Sigma$ for all $d \geq 2$, then
\[ w(\Delta(f),v) = w(\Delta(f)^{(1)},v) + 2 = w(\Delta^{(1)},v) + 2 = w(\Delta,v),\]
where the outer equalities follow from Lemma~\ref{latticewidthproperties}.(iii). \hfill $\blacksquare$\\

A pencil on $C$ that arises as $g_v$ for some lattice direction $v$ is called \emph{combinatorial}.
The number of combinatorial pencils is countable; in fact, by Lemma~\ref{numberofdirections}
there is only a finite number of combinatorial pencils of each given degree. Note that the minimal
degree of a combinatorial pencil is $\text{lw}(\Delta)$, from which we 
immediately find that the gonality $\gamma(C)$ of $C$ is bounded from above by $\text{lw}(\Delta)$.
As we will see in Section~\ref{section_gonality}, equality typically holds.

The correspondence between pairs $\pm v$ of lattice directions and combinatorial
pencils is usually $1$-to-$1$, but there are counterexamples. For
instance, let $\Delta$ be a primitive lattice parallelogram, i.e.\ a polygon of the
form $\conv \{ (0,0), v_1, v_2, v_1 + v_2 \}$ for linearly independent primitive vectors $v_1,v_2 \in \mathbb{Z}^2$.
Then 
\[ w(\Delta, v_1) = w(\Delta, v_2) = \left| \det(v_1,v_2) \right|.\]
Assume that $f$ is supported on the vertices of $\Delta$ only, i.e.\
\[ f = c_{0,0} + c_{1,0}(x,y)^{v_1} + c_{0,1}(x,y)^{v_2} + c_{1,1}(x,y)^{v_1 + v_2}, \]
and that the coefficients $c_{i,j}$ are sufficiently generic.
Then the fiber of $U(f) \rightarrow \mathbb{T}^1 : (x,y) \mapsto (x,y)^{v_1}$ above a point $c \in \mathbb{T}^1 \setminus \{ - c_{0,1} c_{1,1}^{-1} \}$ 
matches with the fiber of $U(f) \rightarrow \mathbb{T}^1 : (x,y) \mapsto (x,y)^{v_2}$ above
\[ - \frac{ c_{0,0} + c  c_{1,0} }{c_{0,1} + c c_{1,1}}. \]
From this it follows that $g_{v_1} = g_{v_2}$.
The same construction works
for the primitive lattice triangle $\Delta = \conv \{ (0,0), v_1, v_2 \}$.\\

\noindent \emph{Example with $v_1 = (3,2)$ and $v_2 = (1,0)$.} The graph below shows the (real affine) zero locus of
$f = 3 + x + x^3y^2 - x^4y^2 \in \mathbb{C}[x,y]$. 
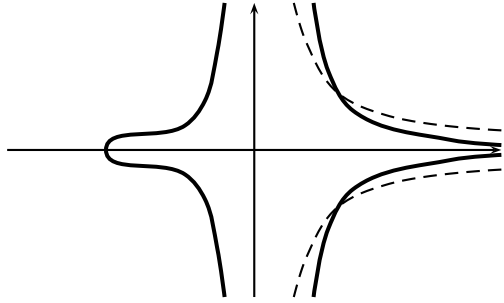
\begin{center}
\psset{unit=0.65cm}
\begin{pspicture}(-5,-3.2)(5,3.2)
\psline{->}(-5,0)(5,0)
\psline{->}(0,-3)(0,3)
\pscurve[linewidth=1.5pt](-0.6,-3)(-0.75,-2)(-1,-1)(-1.4,-0.5)(-3,0)(-1.4,0.5)(-1,1)(-0.75,2)(-0.6,3)
\pscurve[linewidth=1.5pt](5,-0.1)(4,-0.2)(2,-0.8)(1.5,-1.6)(1.2,-3)
\pscurve[linewidth=1.5pt](5,0.1)(4,0.2)(2,0.8)(1.5,1.6)(1.2,3)
\pscurve[linestyle=dashed](5,-0.4)(2.7,-0.7)(1.65,-1.2)(1.05,-2.2)(0.8,-3)
\pscurve[linestyle=dashed](5,0.4)(2.7,0.7)(1.65,1.2)(1.05,2.2)(0.8,3)
\end{pspicture}
\end{center}
The dashed line cuts out a typical fiber of $x^3y^2$, which is also a fiber of $x$.\\

Clearly, by degree considerations, $w(\Delta,v_1) \neq w(\Delta,v_2)$ is a sufficient condition for $v_1,v_2$
to give rise to different combinatorial pencils.
Another sufficient condition is as follows.

\begin{lemma} 
Let $\Delta$ be a two-dimensional lattice polygon
and let $f$
be a $\Delta$-toric Laurent polynomial. Let $v_1 \neq \pm v_2$ be
lattice directions and let $g_{v_1}$ and $g_{v_2}$ be the corresponding combinatorial pencils.
If
\begin{equation} \label{differentpencilssufficientcondition}
  w(\Delta^{(1)}, v_1) \, > \, \left| \det(v_1,v_2) \right| \, - \, 2
\end{equation}
then $g_{v_1} \neq g_{v_2}$.
\end{lemma}

\noindent \textsc{Proof.}
Fibers of
\[ \mathbb{T}^2 \rightarrow \mathbb{T}^1 : (x,y) \mapsto (x,y)^{v_1} \qquad \text{and} \qquad \mathbb{T}^2 \rightarrow \mathbb{T}^1 : (x,y) \mapsto (x,y)^{v_2} \]
intersect each other in at most $\left| \det(v_1,v_2) \right|$ points. Now
because
$\Delta(f)^{(1)} = \Delta^{(1)}$,
condition (\ref{differentpencilssufficientcondition}) implies that
$w(\Delta(f), v_1) > \left| \det(v_1,v_2) \right|$. 
We conclude that
a general fiber of \[ U(f) \rightarrow \mathbb{T}^2 : (x,y) \mapsto (x,y)^{v_1} \]
cannot be contained in a fiber of 
$U(f) \rightarrow \mathbb{T}^2 : (x,y) \mapsto (x,y)^{v_2}$.
The lemma follows.
\hfill $\blacksquare$\\

In the case of lattice width directions we obtain:

\begin{corollary} \label{latticewidthdifferentpencils} Let $\Delta$ be a two-dimensional lattice polygon
and assume that $\sharp (\Delta^{(1)} \cap \mathbb{Z}^2) > 1$. Let $f$
be a $\Delta$-toric Laurent polynomial, let $v_1 \neq \pm v_2$ be
lattice width directions, and let $g_{v_1}$ and $g_{v_2}$ be the according combinatorial pencils.
Then $g_{v_1} \neq g_{v_2}$.
\end{corollary}

\noindent \textsc{Proof.} If $\left| \det(v_1,v_2) \right| = 1$ then
condition (\ref{differentpencilssufficientcondition}) amounts to $\Delta^{(1)} \neq \emptyset$,
which is clearly the case. So by Lemma~\ref{latticewidthproperties}.(iv)
it remains to analyze the case where $\left| \det(v_1,v_2) \right|  = 2$ and $\Delta \cong d\Gamma^5_1$ for
some integer $d \geq 2$ (indeed, $d=1$ is excluded in the statement of the corollary). But here
\[ w(\Delta^{(1)}, v_1) = w((d-1)\Gamma_1^5, v_1) \geq \text{lw}((d-1)\Gamma^5_1) = 2(d-1), \]
so again condition~(\ref{differentpencilssufficientcondition}) is satisfied.  \hfill $\blacksquare$\\

\section{Gonality} \label{section_gonality}

We can now state our refinement of Kawaguchi's theorem \cite[Thm.\,1.3]{Kawaguchi}.

\begin{theorem} \label{maintheorem_gonality}
 Let $f \in k[x^{\pm 1},y^{\pm 1}]$ be non-degenerate with respect to its Newton polygon
 $\Delta = \Delta(f)$. Suppose that $\Delta^{(1)}$ is not equivalent to any of the following:
 \begin{equation} \label{gonality_special_cases}
  \emptyset, \quad (d-3)\Sigma \text{ (for some integer $d \geq 3$)}, \quad \Upsilon, \quad 2\Upsilon, \quad \Gamma_1^5,
  \quad \Gamma_2^5, \quad \Gamma_3^5.
 \end{equation}
 Then every gonality pencil
 on (the smooth projective model of) $U(f)$ is combinatorial.
\end{theorem}

\noindent \emph{Remark.} In case $\Delta^{(1)}$ is among $\Upsilon, 2\Upsilon, \Gamma_1^5,
  \Gamma_2^5, \Gamma_3^5$, there is only a single corresponding $\Delta$,
  namely, $2\Upsilon, 3\Upsilon, 2\Gamma^5_1, 2\Gamma^5_2$ and $2\Gamma^5_3$, respectively.\\

Before we proceed to the proof of Theorem~\ref{maintheorem_gonality}, let us discuss
some corollaries. From the \'enonc\'e it follows that if $\Delta^{(1)}$ is non-equivalent
to any of the polygons listed in (\ref{gonality_special_cases}), then
the gonality of $U(f)$ equals the lattice width of $\Delta$. 
Thus by Lemma~\ref{latticewidthproperties}.(i), if $\Delta^{(1)}$ is not among the polygons listed in (\ref{gonality_special_cases}) then
the gonality of $U(f)$ equals $\text{lw}(\Delta^{(1)}) + 2$. The other instances of $\Delta^{(1)}$ can be analyzed case by case:
\begin{itemize}
  \item If $\Delta^{(1)} = \emptyset$ then $U(f)$ is rational, hence of gonality $1$.
  \item If $\Delta^{(1)} \cong (d-3)\Sigma$ then $U(f)$ is birationally equivalent to a smooth projective plane curve of degree $d$, hence of gonality $d-1$ by a result of Namba \cite{Namba}
        (a proof can also be found in \cite[Prop.\,3.13]{Serrano}).
  \item If $\Delta^{(1)} \cong \Upsilon$ then $U(f)$ is a non-hyperelliptic curve of genus $4$, hence of
        gonality $3$.
  \item If $\Delta^{(1)} \cong 2\Upsilon$ then $U(f)$ is birationally equivalent to a smooth intersection of two cubics in $\mathbb{P}^3$, hence of
  gonality $6$ by a result of Martens (see \cite[Thm.\,9]{CaCo} for more details).
  \item If $\Delta^{(1)} \cong \Gamma^5_i$ ($i=1,2,3$) then $U(f)$ is a non-hyperelliptic, non-trigonal curve
  of genus $5$ by \cite[Lem.\,3]{CaCo}, hence of gonality $4$.
\end{itemize}

We conclude:
\begin{corollary} \label{cacoconjecture}
 Let $f \in k[x^{\pm 1},y^{\pm 1}]$ be non-degenerate with respect to its Newton polygon
 $\Delta = \Delta(f)$. Then the gonality of $U(f)$ equals $\emph{lw}(\Delta^{(1)}) + 2$,
 unless $\Delta^{(1)} \cong \Upsilon$ (i.e.\ $\Delta \cong 2\Upsilon$), in which case it equals $3$.
\end{corollary}

Unless $\Delta^{(1)} \cong \Upsilon$ we can even read off the number of gonality pencils:

\begin{corollary} \label{numberofpencils}
 Let $f \in k[x^{\pm 1}, y^{\pm 1}]$ be non-degenerate with respect to its Newton polygon
 $\Delta = \Delta(f)$. 
 \begin{itemize}
  \item If $\Delta^{(1)} = \emptyset$ then there is a unique gonality pencil.
  \item If $\Delta^{(1)} \cong \Upsilon$ then the number of gonality pencils is at most $2$.
  \item If $\Delta^{(1)} \cong (d-3)\Sigma$ for some $d \geq 3$, or if
  $\Delta^{(1)} \cong 2\Upsilon, \Gamma_1^5, \Gamma_2^5, \Gamma_3^5$, then there
  are infinitely many gonality pencils.
  \item In all other cases the number of gonality pencils equals the number of lattice width directions.
  In particular, the number of gonality pencils is at most $4$, and the bound is met iff
  $\Delta^{(1)} \cong d\Gamma_1^5$ for some $d \geq 2$.
 \end{itemize}
\end{corollary}

\noindent \textsc{Proof.} The first three claims follow
from the considerations above: rational curves have a unique gonality pencil, 
non-hyperelliptic genus $4$ curves carry one or two $g^1_3$'s \cite[Ex.\,IV.5.5.2]{hartshorne},
smooth plane degree $d$ curves admit infinitely many $g^1_{d-1}$'s \cite[Prop.\,3.13]{Serrano},
smooth intersections of cubics in $\mathbb{P}^3$ carry infinitely many $g^1_6$'s \cite[pp.\,174-175]{ELMS},
and non-hyperelliptic, non-trigonal curves of genus $5$ have infinitely many $g^1_4$'s \cite[Ex.\,IV.F]{Arbarello}. The last
claim follows from Theorem~\ref{maintheorem_gonality}, combined with Lemma~\ref{latticewidthproperties}.(iv) and 
Corollary~\ref{latticewidthdifferentpencils}.
\hfill $\blacksquare$\\

\noindent \emph{Example (revisited, see Section~\ref{section_combpencils}).} Let $\Delta$ be the lattice polygon
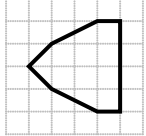
\begin{center}
\psset{unit=0.3cm}
\begin{pspicture}(-1,-3)(5,3)
\psgrid[subgriddiv=1,griddots=10,gridcolor=lightgray,gridlabels=0pt](-1,-3)(5,3)
\pspolygon[linewidth=1.5pt](0,0)(1,-1)(3,-2)(4,-2)(4,2)(3,2)(1,1)
\end{pspicture}
\end{center}
and let $f \in k[x^{\pm 1}, y^{\pm 1}]$ be a $\Delta$-non-degenerate (or $\Delta$-toric) 
Laurent polynomial. Then $U(f)$ is a $4$-gonal genus $7$ curve carrying exactly
two $g^1_4$'s.\\

\noindent \emph{Remarks.}
\begin{itemize}
\item Corollary~\ref{cacoconjecture} implies a conjecture by the current
authors \cite[Conj.~1]{CaCo}. It does not imply the corresponding conjecture
on metric graphs \cite[Conj.~3 + Err.]{CaCo}.
\item Corollary~\ref{cacoconjecture} also implies that if $\Delta^{(1)} \cong \Upsilon$ (i.e.\ if $\Delta \cong 2\Upsilon$), then a combinatorial gonality
pencil cannot exist. The same conclusion holds for $\Delta \cong d\Sigma$ for $d \geq 2$. In all other
cases, there exists at least one combinatorial gonality pencil.
\item In case $\Delta^{(1)} \cong \Upsilon$ then both one and two $g^1_3$'s can occur,
depending on whether the quadric on which the curve canonically embeds is singular or not \cite[Ex.\,IV.5.5.2]{hartshorne};
see \cite[Thm.\,4]{CaCoCanonical} for an explicit description of this quadric.
\item Let $k'$ be an arbitrary field of characteristic $0$ with algebraic closure
$k$. Suppose that
$f \in k'[x^{\pm 1}, y^{\pm 1}]$ is non-degenerate with respect to its Newton polygon
when considered
as a Laurent polynomial over $k$. If $\Delta(f) \not \cong 2\Upsilon, d\Sigma$ then
the above remark implies
that $\gamma(U(f)) = \gamma_{k'}(U(f))$, where $\gamma_{k'}(U(f))$ is the minimal
degree of a \emph{$k'$-rational} map to $\mathbb{P}^1$. If $\Delta(f) \cong 2\Upsilon$ or $\Delta(f) \cong d\Sigma$ for some $d \geq 2$ then this
may not be true. (Example: $x^2 + y^2 + 1 \in \mathbb{R}[x,y]$.)
\item By letting $k' = \mathbb{C}(\!(t)\!)$, the preceding remark lends prudent support in 
favor of a conjecture by M.\,Baker (stating that
the gonality of a graph equals the gonality of the associated metric graph \cite[Conj.~3.14]{baker}) in
the case of graphs associated to regular subdivisions of lattice polygons \cite[Err.~\S1]{CaCo}.
\item If $\Delta^{(1)}$ is neither among
the polygons excluded in Theorem~\ref{maintheorem_gonality}, nor of the form
$d\Gamma_1^5$ for some $d \geq 2$, then 
Lemma~\ref{latticewidthproperties}.(v) 
implies that two different gonality pencils on $U(f)$ are always \emph{independent}, in the sense
that they span a base-point free linear system of rank $2$, defining a morphism
$U(f) \rightarrow \mathbb{P}^2$ that induces a birational equivalence between $U(f)$ and its image.
(For general lattice directions $v_1 \neq \pm v_2$ the morphism $C \rightarrow \mathbb{P}^2$ 
defined by $g_{v_1}$ and $g_{v_2}$ induces a degree $\left| \det(v_1,v_2) \right| $ cover.)
See \cite[(1.2)]{Coppens} for more background on this terminology.\\
\end{itemize}
\vspace{-0.3cm}

We now give a proof of Theorem~\ref{maintheorem_gonality}. We recall that the main ideas are
taken from Kawaguchi \cite{Kawaguchi}, but that our proof covers the case where $U(f)$
is birationally equivalent to a smooth projective plane curve (the key ingredient here being
the block of text surrounding (\ref{bezout}) below).\\

\noindent \textsc{Proof of Theorem~\ref{maintheorem_gonality}.} Let $g = \sharp(\Delta^{(1)} \cap \mathbb{Z}^2)$ be the geometric genus of $U(f)$.
Note that our assumptions imply $g \geq 2$. Recall that $\text{lw}(\Delta^{(1)}) = 0$ if and only if $U(f)$ is hyperelliptic.
By Lemma~\ref{latticewidthproperties}.(i) this holds if and only if $\text{lw}(\Delta) = 2$, hence
a $g^1_2$ can be computed by projection along some lattice direction. Since the $g^1_2$ of a hyperelliptic curve
is unique, Theorem~\ref{maintheorem_gonality} follows in this case. Thus we may assume that $\Delta^{(1)}$ is two-dimensional
and that $U(f)$ is of gonality $\gamma \geq 3$.

From Lemma~\ref{cruciallemma} we know that $f$ is $\Delta^\text{max}$-toric,
i.e.\ $\varphi_{\Delta^\text{max}}(U(f))$ completes to a smooth
projective curve $C \subset \text{Tor}(\Delta^\text{max})$. Let $\Sigma'$ be a minimal smooth subdivision
of $\Sigma_{\Delta^\text{max}}$ and let $\mu : \text{Tor}(\Sigma') \rightarrow \text{Tor}(\Sigma_{\Delta^\text{max}})$ be
the corresponding birational morphism. Let $C'$ be the strict transform of $C$ under $\mu$. 
Because the smooth subdivision was chosen minimal,
$C'$ does not meet the exceptional locus of $\mu$. In particular, $\mu|_{C'}$ is an isomorphism of curves
and $\Delta_{C'} \cong_t \Delta^\text{max}$. Since $\text{Tor}(\Sigma')$ is smooth, every Weil divisor is Cartier.

By the BKK theorem (recall that $C'$ is a convex divisor),
\[ C'^2  =  \text{MV}(\Delta^\text{max}, \Delta^\text{max}) = 2 \text{Vol}(\Delta^\text{max})
   \geq \frac{3}{4} \text{lw}(\Delta^\text{max})^2 = \frac{3}{4} \text{lw}(\Delta)^2,\]
where the third and fourth (in)equalities follow from Lemma~\ref{latticewidthproperties}.(i,vi).
For small lattice widths this bound can be improved: using the
data from \cite{movingout} one can computationally verify that
\begin{equation} \label{computationalbounds}
 C'^2 = 2 \text{Vol}(\Delta^\text{max}) \geq \left\{ \begin{array}{ll} 18 & \text{if $\text{lw}(\Delta^\text{max}) = 3$}, \\
                           20 & \text{if $\text{lw}(\Delta^\text{max}) = 4$}, \\
                           25 & \text{if $\text{lw}(\Delta^\text{max}) = 5$}, \\
                           28 & \text{if $\text{lw}(\Delta^\text{max}) = 6$} \\
         \end{array}
 \right.
\end{equation}
(remark that by Pick's theorem it suffices to verify these inequalities for small genus only).
Magma code assisting the reader in this can be found in the accompanying file \verb"gonal.m". 
The patient reader can also do an elaborate analysis by hand, following Kawaguchi \cite[Props.~3.10--3.12,4.3]{Kawaguchi}.
We stress that for these bounds it is essential that $\Delta^\text{max}$ is maximal 
and that $\Delta^{(1)}$ is not among the polygons listed in (\ref{gonality_special_cases}).

We now come to the heart of the proof. Fix a gonality pencil $g^1_\gamma$ and let
$p : C' \rightarrow \mathbb{P}^1$ be a corresponding morphism of degree $\gamma$.
A theorem by Serrano \cite[Thm.\,3.1]{Serrano} states that if
$C'^2 > (\gamma + 1)^2$ then $p$ can be extended to a morphism $\text{Tor}(\Sigma') \rightarrow \mathbb{P}^1$.
From this it will follow that $p$ is combinatorial (as explained in the last paragraph of the proof). Unfortunately, in general we only have 
that\footnote{But note that for `most' lattice polygons, the stronger bound $C'^2 > (\gamma + 1)^2$ \emph{does} hold, in which case the proof
simplifies a lot.}
\begin{equation} \label{selfintlowbound}
 C'^2 \geq \frac{3}{4} \text{lw}(\Delta)^2 \geq \frac{3}{4} \gamma^2.
\end{equation}
To bridge this, we follow an approach of Harui \cite{Harui}, who dug into Serrano's proof to extract
Theorem~\ref{thm_serrano} below.

We proceed by contradiction: assume that $p$ cannot be extended to all of $\text{Tor}(\Sigma')$. Then
by Theorem~\ref{thm_serrano} (note that $C'^2 > 4\gamma$)
there exists an effective divisor $V$ on $\text{Tor}(\Sigma')$ satisfying
\begin{equation} \label{selfintuppbound}
 1 \leq s < C' \cdot V - s \leq \gamma \qquad \text{and} \qquad C'^2 \leq \frac{(\gamma + s)^2}{s},
\end{equation}
where $s = V^2$. We may assume that $V$ is torus-invariant, i.e.\ $V = \sum a_\ell D_\ell$ for certain
integers $a_\ell$ (where the $D_\ell$'s are the torus-invariant prime divisors of $\text{Tor}(\Sigma')$). From our bounds (\ref{selfintlowbound}) and (\ref{selfintuppbound})
we see that 
\[ \frac{3}{4} \gamma^2 \leq \frac{(\gamma + s)^2}{s} < \frac{(\gamma + \gamma)^2}{s} = \frac{4}{s} \gamma^2 \] which implies that $s \leq 5$. Rewrite the first inequality as
$(3s - 4)\gamma^2 - 8s \gamma - 4s^2 \leq 0$: if $s \geq 2$ then the largest real root of the left-hand side, when viewed
as a polynomial in $\gamma$, is given by
$(4s + 2s \sqrt{3s})/(3s - 4)$ which for $s \leq 5$ is strictly less than $9$. 
We conclude that if $\gamma \geq 9$ then $s=1$. A finer analysis
using the better bounds (\ref{computationalbounds}) shows
that $\gamma \geq 4$, and that $s = 1$ except possibly if $\gamma \in \{6,7,8\}$ in which case $s \in \{1,2\}$.

We claim that this implies $h^0(\text{Tor}(\Sigma'),V) \leq s+1$. Suppose not, then $\Delta_V$ contains at least $s+2$ lattice points.
  Let $\Gamma$ be the convex hull of the lattice points in $\Delta_V$ and let $D_\Gamma = \sum_\ell a'_\ell D_\ell$ be the
  torus-invariant divisor obtained by taking the $a'_\ell$'s minimal such that
  \[ \Gamma \subset \left\{ \, \left. (i,j) \in \mathbb{R}^2 \, \right| \, <\! (i,j), v_\ell \! > \ \geq - a'_\ell \, \right\}.\]
  One verifies
  that $D_\Gamma$ is convex, that $\Delta_{D_\Gamma} = \Gamma$, and that $a_\ell \geq a'_\ell$ for all $\ell$, i.e.\ $V - D_\Gamma$ is effective.
  \begin{itemize}
    \item Suppose that, up to a unimodular transformation, $\Gamma$ contains a horizontal line segment $I$ of
    length $\geq 2$. Then $C' \cdot V = C' \cdot (D_\Gamma + (V - D_\Gamma))$ is bounded from below by
    \[ C' \cdot D_\Gamma = \text{MV}(\Delta^\text{max},\Gamma) \geq \text{MV}(\Delta^\text{max},I) \geq 2 \, \text{lw}(\Delta^\text{max}) = 2 \, \text{lw}(\Delta) \geq 2 \gamma, \]
    where the first inequality follows because $\text{MV}$ is an increasing function.This contradicts $C' \cdot V \leq \gamma + s$.
    \item So we can assume that $\Gamma$ does not contain such a line segment.
    \begin{itemize}
    \item
    If $s=1$ then $\Gamma$ contains at least $3$ non-collinear lattice points. But 
    then, by applying a unimodular transformation if needed,
    we may assume that $\Sigma \subset \Gamma$. One finds
  \begin{equation} \label{bezout}
   \qquad \qquad \quad \quad C' \cdot V \geq \text{MV}(\Delta^\text{max},\Gamma) \geq \text{MV}(\Delta^\text{max},\Sigma) = d
  \end{equation}
  where $d$ is the smallest integer such that
  $\Delta^\text{max}$ is contained in a translate of $d\Sigma$ (indeed, this follows from B\'ezout's theorem). Then
  $\Delta^{(1)} \subset (d-3)\Sigma$, and by our assumptions this inclusion is strict. It follows
  that $\text{lw}(\Delta^{(1)}) \leq d-4$, hence by Lemma~\ref{latticewidthproperties}.(i) that $\text{lw}(\Delta) = \text{lw}(\Delta^\text{max}) \leq d-2$.
  From (\ref{bezout}) we conclude that $C' \cdot V \geq \text{lw}(\Delta) + 2$. This contradicts
  $C' \cdot V \leq \gamma + 1$.
  \item If $s=2$ then $\Gamma$ contains at least $4$ lattice points. By our assumption that it contains no line
  segment of integral length $2$, we can assume $\square \subset \Gamma$ or $\Upsilon \subset \Gamma$, again by applying a unimodular transformation if needed.
  In the former case we have
  \[ \qquad \qquad \quad \quad C' \cdot V \geq \text{MV}(\Delta^\text{max},\Gamma) \geq \text{MV}(\Delta^\text{max},\square) = a + b, \]
  where $(a,b)$ is the `bidegree' of $f$, i.e.\ the minimal couple of values for which $\Delta^\text{max}$
  is contained in a translate of $[0, a] \times [0, b]$. This follows
  from the BKK theorem applied to $\text{Tor}(\square) = \mathbb{P}^1 \times \mathbb{P}^1$, and implies that $C' \cdot V \geq 2 \, \text{lw}(\Delta) \geq 2 \gamma$.
  In the latter case, by the BKK theorem applied to $\text{Tor}(\Upsilon)$, one similarly finds
    \[  \qquad \qquad \quad \quad C' \cdot V \geq \text{MV}(\Delta^\text{max},\Gamma) \geq \text{MV}(\Delta^\text{max},\Upsilon) = 3d, \]
   where $d$ is the smallest integer such that $\Delta^\text{max}$ is contained in a translate of $d\Upsilon$.
  Because 
  \begin{equation} \label{willbesharpened}
    \qquad \qquad \quad \quad 2d = \text{lw}(d \Upsilon) \geq \text{lw}(\Delta^\text{max}) \geq \gamma 
  \end{equation}  
  we find that $C' \cdot V \geq \frac{3}{2}\gamma$.
  In both cases this contradicts $C' \cdot V \leq \gamma + 2$; recall that $\gamma \geq 6$ in the $s=2$ case. \vspace{-0.2cm}\\
  
  \noindent \emph{Remark.} The bound $C' \cdot V \geq \frac{3}{2}\gamma$ can be proven more easily by noting
  that $\Upsilon$ contains a line segment of integral length $\frac{3}{2}$; however, 
  the argument using (\ref{willbesharpened}) will reappear in the proof of Theorem~\ref{maintheorem_neargonal},
  so we have included it for the sake of consistency.\\
  \vspace{-0.2cm}
  
  \end{itemize}
  \end{itemize}
  Our claim that $h^0(\text{Tor}(\Sigma'),V) \leq s+1 $ follows.

Because a lattice
polygon having at most $3$ lattice points cannot have any interior lattice points,
we deduce that $h^0(\text{Tor}(\Sigma'), V + K) = 0$, with $K =  - \sum_\ell D_\ell$ the canonical
divisor from Section~\ref{section_toric}. The Riemann-Roch theorem yields that
\[ \frac{1}{2}(V+K) \cdot V = \qquad \qquad \qquad \qquad \qquad \qquad \qquad \qquad \qquad \qquad \qquad \qquad \qquad \]
\[h^0(\text{Tor}(\Sigma'),V+K) - h^1(\text{Tor}(\Sigma'),V+K) + h^0(\text{Tor}(\Sigma'),-V) - \chi(\mathcal{O}_\text{Tor}(\Sigma'))\]
is bounded by $-\chi(\mathcal{O}_\text{Tor}(\Sigma')) = -1$, i.e.\ $K \cdot V \leq -s-2$. But then Riemann-Roch
also tells us that
\[
h^0(\text{Tor}(\Sigma'),V)  =   h^1(\text{Tor}(\Sigma'),V) - h^0(\text{Tor}(\Sigma'),K-V) + \frac{1}{2}V \cdot (V - K) + 1\]
is at least $s + 2$.

Thus we run into the desired contradiction, and we conclude that
$p$ can be extended to all of $\text{Tor}(\Sigma')$.
Let $\tilde{p} : \text{Tor}(\Sigma') \rightarrow \mathbb{P}^1$ be such that $\tilde{p}|_{C'} = p$.
Let $F$ be a fiber of $\tilde{p}$, so that $F \cdot C' = \gamma$. Then $C' \cdot (F - C') \leq \gamma - \frac{3}{4} \gamma^2 < 0$.
Since $C'$ is nef it follows that $h^0(\text{Tor}(\Sigma'),F-C') = 0$.
Now by tensoring
the short exact sequence
\[ 0 \rightarrow \mathcal{O}_{\text{Tor}(\Sigma')}(- C') \rightarrow \mathcal{O}_{\text{Tor}(\Sigma')} \rightarrow \mathcal{O}_{C'} \rightarrow 0 \]
with $\mathcal{O}_{\text{Tor}(\Sigma')}(F)$ and taking cohomology,
we find the exact sequence
\[ 0 \rightarrow H^0(\text{Tor}(\Sigma'), F - C') \rightarrow H^0(\text{Tor}(\Sigma'),F) \rightarrow H^0(C',F|_{C'}) \rightarrow \dots, \]
which proves that $h^0(\text{Tor}(\Sigma'), F) \leq 2$; here we used that $h^0(C',F|_{C'}) = 2$ because $g^1_\gamma$ is complete.
Thus $|F|$ is a linear system of rank $1$, i.e.\ every element of $|F|$ is a fiber of $\tilde{p}$.
Let $D$ be a torus-invariant divisor that is equivalent to $F$. By translating if necessary
we may assume that $(0,0) \in \Delta_D$, so that $D$ is effective. But then $D \in |F|$ and
\[ H^0(\text{Tor}(\Sigma'),D) = \left\langle 1, x^ay^b \right\rangle \]
for some primitive $(a,b) \in \mathbb{Z}^2$. We find
that $\tilde{p}|_{\mathbb{T}^2} : (x, y) \mapsto x^a y^b$ (up to an automorphism of $\mathbb{P}^1$), i.e.\
$g^1_\gamma = g_{(a,b)}$.
 \hfill $\blacksquare$\\

\begin{theorem}[Serrano, 1987] \label{thm_serrano}
Let $C$ be a smooth projective curve on a smooth projective surface $S$, and let
$p : C \rightarrow \mathbb{P}^1$ be a surjective morphism of degree $d$. Suppose
that $C^2 > 4d$ and that $p$ cannot be extended to a morphism $S \rightarrow \mathbb{P}^1$.
Then there exists an effective divisor $V$ on $S$ for which
\[ 0 < V^2 < V \cdot (C - V) \leq d \qquad \text{and} \qquad C^2 \leq \frac{(d + V^2)^2}{V^2}.\]
\end{theorem}

\noindent \textsc{Proof.} By contradiction. Suppose that such an effective divisor $V$ does not exist,
then one can replace \emph{Claim 6} in Serrano's proof \cite[p.\ 401]{Serrano} by the following reasoning
(the text below does not make sense without Serrano's paper at hand):
\begin{adjustwidth}{2.5em}{0pt}
   \emph{Claim 6: $a = 0$.} Suppose that $a > 0$. Then $V_1$ is an effective divisor
   such that $0 < V_1^2 < V_1 \cdot V_2 \leq d$ because $a < e$.
   On the other hand,
   \[ C^2 = a + 2e + b  \leq a + 2e + \frac{e^2}{a} \leq a + 2d + \frac{d^2}{a} = \frac{(a + V_1^2)^2}{V_1^2}.\]
   Since $V_2 = C - V_1$ this contradicts our hypothesis. Hence, $a = 0$.
\end{adjustwidth}
The rest of the proof can be copied word by word. \hfill $\blacksquare$

\section{Near-gonal pencils} \label{section_neargonal}

By a \emph{near-gonal pencil} on a smooth projective curve $C/k$ we mean a base-point free $g^1_{\gamma(C) + 1}$ (note 
that such pencils need not exist). 
The method of the previous section can be adapted to show that, apart from some reasonably well-understood exceptional instances
of $\Delta$, every near-gonal pencil on a $\Delta$-non-degenerate curve is
combinatorial.

It is convenient to state our main result in terms of the \emph{lattice size},
a notion to which we have devoted a separate paper \cite{latticesize}.
If $\Delta \neq \emptyset$, then its lattice size is defined as the minimal integer $d \geq 0$
such that $\Delta$ is equivalent to a lattice polygon
that is contained in $d\Sigma$. We denote this integer by $\lsize(\Delta)$, and let $\lsize(\emptyset) = -2$.
If $\Delta$ is two-dimensional then, as in the case of 
the lattice width (cf.\ Lemma~\ref{latticewidthproperties}.(i)), there exists
an expression for $\lsize(\Delta)$ in terms of $\lsize(\Delta^{(1)})$, allowing one
to compute $\lsize(\Delta)$ by gradually peeling off the polygon \cite[Thm.\,3.5]{latticesize}. For our needs,
one of the main results of \cite{latticesize}
can be reformulated as follows:
\begin{theorem} \label{s2upper}
 Let $f \in k[x^{\pm 1}, y^{\pm 1}]$ be non-degenerate
 with respect to its Newton polygon $\Delta = \Delta(f)$. Then the
 minimal degree of a (possibly singular) projective plane curve
 that is birationally equivalent to $U(f)$ is bounded by
 $\lsize(\Delta^{(1)}) + 3$.
 If $\Delta^{(1)} \cong (d-1)\Upsilon$ for a certain integer $d \geq 2$ (i.e.\ $\Delta \cong d\Upsilon$), 
 then it is moreover bounded by $3d - 1$.
\end{theorem}
\noindent \textsc{Proof.} See \cite[Thm.\,1.3]{latticesize}. \hfill $\blacksquare$\\

\noindent \emph{Remarks.} 
\begin{itemize}
 \item If $\Delta^{(1)} \cong (d-1)\Upsilon$ then $\lsize(\Delta^{(1)}) + 3 = 3d$ (as can be
 verified using \cite[Thm.\,3.5]{latticesize}). So the second bound is sharper in this case.
 \item We expect that the (smallest applicable) bound of Theorem~\ref{s2upper}
is in fact sharp; see \cite[\S7]{latticesize} for a discussion.\\
\end{itemize}

Our main result is as follows:
\begin{theorem} \label{maintheorem_neargonal}
 Let $f \in k[x^{\pm 1}, y^{\pm 1}]$ be non-degenerate with respect to its Newton polygon $\Delta = \Delta(f)$, 
 and let $\gamma$ be the gonality of $U(f)$.
 Suppose 
 that
 \begin{equation} \label{planecurvecondition}
  \lsize(\Delta^{(1)}) \geq \emph{lw}(\Delta^{(1)}) + 2
 \end{equation}
 and that $\Delta^{(1)} \not \cong 2\Upsilon, 3\Upsilon, \Gamma^7, \Gamma^8$.
 Then every base-point free $g^1_{\gamma + 1}$ on the smooth projective model of
 $U(f)$ is combinatorial.
\end{theorem}

Before we proceed to the proof of Theorem~\ref{maintheorem_neargonal}, let us 
analyze the list of polygons that are excluded in the statement
(this is a strict extension of the list of polygons that were exluded in the statement of 
Theorem~\ref{maintheorem_gonality}).
First note that Theorem~\ref{maintheorem_neargonal} implies
that if $\Delta$ is not among the excluded polygons, the number
of base-point free $g^1_{\gamma + 1}$'s is finite. Opposed to that, we have:

\begin{lemma}
 If $\Delta$ violates condition (\ref{planecurvecondition}) or $\Delta^{(1)} \cong 2\Upsilon, \Gamma^7$,
 then the number of base-point free $g^1_{\gamma + 1}$'s is infinite.
\end{lemma}

\noindent \textsc{Proof.}
A violation of condition (\ref{planecurvecondition}) implies that $U(f)$ is birationally
equivalent to a (possibly singular) plane curve of degree at most $\gamma + 2$. Indeed:
\begin{itemize}
  \item If $\Delta^{(1)} \cong \Upsilon$ 
 then $U(f)$ is a non-hyperelliptic genus $4$ curve, hence of gonality $3$. It is 
 known that such curves admit a plane model of degree $5$; see e.g.\ \cite[Ex.\,IV.5.4]{hartshorne}. 
 \item If $\Delta^{(1)} \not \cong \Upsilon$ but $\lsize(\Delta^{(1)}) < \text{lw}(\Delta^{(1)}) + 2$, 
 then by Corollary~\ref{cacoconjecture} the assumption 
 can be rephrased as $\lsize(\Delta^{(1)}) < \gamma$. Along
 with Theorem~\ref{s2upper} this implies that
 $U(f)$ has a projective plane model of degree at most $\gamma + 2$.
\end{itemize}
It follows that $U(f)$ must have a plane model of degree exactly $\gamma + 1$ or $\gamma + 2$, because
a model of degree at most $\gamma$ would contradict that $\gamma$ equals the gonality (by projecting from
a point on this plane model). But then
there exist infinitely many base-point free $g^1_{\gamma+1}$'s, obtained either by projection from
a point outside the plane model, or by projection from a non-singular point on the plane model.

If $\Delta^{(1)} \cong 2\Upsilon$, so that $\Delta \cong 3\Upsilon$,
then $U(f)$ is a $6$-gonal curve that is birationally
 equivalent to a smooth intersection of two cubics 
 in $\mathbb{P}^3 = \text{Proj} \, k[X_{0,0}, X_{-1,-1}, X_{1,0}, X_{0,1}]$,
where one of the cubics is just $\text{Tor}(\Upsilon)$, i.e.\ it is given by $X_{0,0}^3 - X_{-1,-1} X_{1,0} X_{0,1}$
 (see the according
 remark following Theorem~\ref{maintheorem_gonality}).
 By the trisecant lemma \cite[IV.Prop.\,3.8 and IV.Thm.\,3.9]{hartshorne} we can
 find a point on this curve, the general secant through which is not a multisecant.
 Projecting from this point gives a birational equivalence with a plane curve of degree $8$, and
 hence we again obtain infinitely many $g^1_7$'s.
 
 Finally, if $\Delta^{(1)} \cong \Gamma^7$ then $\gamma = 4$. Now there exists at least
 one base-point free $g^1_5$ (namely $g_{(0,1)}$). By
Brill-Noether theory it then follows that the number of base-point free $g^1_5$'s is infinite. 
\hfill $\blacksquare$\\

The exclusion of $3\Upsilon$ (in which case $\gamma = 8$) is also necessary:

\begin{lemma}
 If $\Delta^{(1)} \cong 3\Upsilon$ then there exists a base-point free
 $g^1_9$, while
 there are no combinatorial $g^1_9$'s.
\end{lemma}

\noindent \textsc{Proof.} 
If $\Delta^{(1)} \cong 3\Upsilon$, then $\Delta \cong 4\Upsilon$ and $U(f)$ is a curve of genus $19$ that is
birationally
equivalent to a smooth intersection of $\text{Tor}(\Upsilon)$ and a quartic  
in $\mathbb{P}^3$.
By Theorem~\ref{maintheorem_gonality} our curve is $8$-gonal, and there are exactly three $g^1_8$'s. 
Geometrically, the three $g^1_8$'s can be visualized as pencils of planes through the three lines of $\text{Tor}(\Upsilon)$. 
By the trisecant lemma we can find a point on the
curve that is
\begin{enumerate}
  \item not contained in any of these three lines, and
  \item the general secant line through which is not a multisecant.
\end{enumerate}
Projecting from such a point gives a birational map to a 
plane curve of degree $11$, the map being birational because of condition~(2).
Genus considerations yield that the curve must be singular. Moreover, the singular points all 
have multiplicity $2$. Indeed, if there were
a singularity of multiplicity $3$, the pencil of lines through this point would
cut out one of our $g^1_8$'s, which is impossible by condition~(1). On the other hand, a singularity of
higher multiplicity would contradict that the gonality is $8$. Then projecting from
such a singular point of multiplicity $2$ yields a base-point free $g^1_9$.
We leave it to the reader to verify that there are indeed no combinatorial $g^1_9$'s. \hfill $\blacksquare$\\

Finally, if $\Delta^{(1)} \cong \Gamma^8$, so that
$\Delta \cong \conv \{ (0,0), (6,2), (2,4) \}$, then $\gamma = 4$ and it can
be checked that
there are no combinatorial $g^1_5$'s. On the other hand, 
the Laurent polynomial
$f = 1 - x^6y^2 - x^2y^4$ is non-degenerate with respect to its Newton polygon, 
while $U(f)$ admits a rational map
\[ U(f) \rightarrow \mathbb{A}^1 : (x,y) \mapsto \frac{1 - xy^2}{x^3y} \]
of degree $5$, and therefore carries
a base-point free $g^1_5$. Moduli-theoretic considerations then allow one
to draw the same conclusion for a non-empty open subset of the space of Laurent polynomials $f \in k[x^{\pm 1}, y^{\pm 1}]$
that are supported on $\Delta$. Unfortunately, this does not prove the corresponding statement for all $\Delta$-non-degenerate (or
$\Delta$-toric) Laurent polynomials, even though we believe that it should be true. But
in any case this shows that the exclusion of $\Gamma^8$ is also necessary.


We now prove Theorem~\ref{maintheorem_neargonal}:\\

\noindent \textsc{Proof of Theorem~\ref{maintheorem_neargonal}.} 
This is very similar to the proof of Theorem~\ref{maintheorem_gonality}.
Let $g = \sharp(\Delta^{(1)} \cap \mathbb{Z}^2)$ be the geometric
genus of $U(f)$. The assumptions imply that $g \geq 3$. 
Because hyperelliptic curves of genus at least $3$ never carry
a base-point free $g^1_3$, we can assume
that $\Delta^{(1)}$ is two-dimensional and that $U(f)$ is of gonality 
$\gamma \geq 3$.

As before, let $C$ be the Zariski closure of $\varphi_{\Delta^\text{max}}(U(f))$ inside
$\text{Tor}(\Delta^\text{max})$, let $\Sigma'$ be a minimal smooth 
subdivision of $\Sigma_{\Delta^\text{max}}$, let $\mu : \text{Tor}(\Sigma') \rightarrow \text{Tor}(\Sigma_{\Delta^\text{max}})$ be
the corresponding birational morphism, and let $C'$ be the strict transform of $C$ under $\mu$.
Recall that $C'^2 \geq \frac{3}{4}\gamma^2$.
Using the data from \cite{movingout}, our list of sharpened lower bounds (\ref{computationalbounds})
can be adapted and extended to
\begin{equation} \label{computationalbounds2}
 C'^2 = 2 \text{Vol}(\Delta^\text{max}) \geq \left\{ \begin{array}{ll} 
                           24 & \text{if $\text{lw}(\Delta^\text{max}) = 3$}, \\                    
                           24 & \text{if $\text{lw}(\Delta^\text{max}) = 4$}, \\                    
                           30 & \text{if $\text{lw}(\Delta^\text{max}) = 5$}, \\                    
                           34 & \text{if $\text{lw}(\Delta^\text{max}) = 6$}, \\                    
                           46 & \text{if $\text{lw}(\Delta^\text{max}) = 7$}, \\
                           55 & \text{if $\text{lw}(\Delta^\text{max}) = 8$}, \\
         \end{array}
 \right.
\end{equation}
unless $\Delta^\text{max}$ is equivalent to one of the following three polygons,
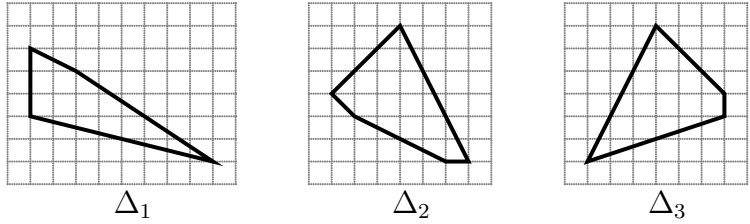
\begin{center}
\psset{unit=0.3cm}
\begin{pspicture}(-3,-4)(7,5.4)
\psgrid[subgriddiv=1,griddots=10,gridcolor=gray,gridlabels=0pt](-3,-3)(7,5)
\pspolygon[linewidth=1.5pt](-2,0)(6,-2)(0,2)(-2,3)
\rput(2.5,-3.9){$\Delta_1$}
\end{pspicture}
\qquad
\begin{pspicture}(-3,-4)(5,5.4)
\psgrid[subgriddiv=1,griddots=10,gridcolor=gray,gridlabels=0pt](-3,-3)(5,5)
\pspolygon[linewidth=1.5pt](-2,1)(-1,0)(3,-2)(4,-2)(1,4)
\rput(1.5,-3.9){$\Delta_2$}
\end{pspicture}
\qquad
\begin{pspicture}(-3,-4)(5,5.4)
\psgrid[subgriddiv=1,griddots=10,gridcolor=gray,gridlabels=0pt](-3,-3)(5,5)
\pspolygon[linewidth=1.5pt](-2,-2)(4,0)(4,1)(1,4)
\rput(1.5,-3.9){$\Delta_3$}
\end{pspicture}
\end{center}
whose respective lattice widths and doubled volumes are $5,6,6$ and $25,32,33$.
See the accompanying Magma file \verb"neargonal.m" for assistance
in verifying these bounds.
It is again essential that $\Delta^\text{max}$ is maximal
and that $\Delta^{(1)}$ is not among the polygons excluded in the \'enonc\'e (recall 
that this is a strict extension of the list of polygons that were excluded in Theorem~\ref{maintheorem_gonality}).

For now, assume that $\Delta^\text{max} \not \cong \Delta_1, \Delta_2, \Delta_3$: we will deal with
these polygons later.
Consider a base-point free $g^1_{\gamma + 1}$ on $C'$ and 
let $p : C' \rightarrow \mathbb{P}^1$ be a corresponding morphism of degree $\gamma + 1$ (which exists
precisely because our $g^1_{\gamma + 1}$ is base-point free).
Assume
that $p$ cannot be extended to all of $\text{Tor}(\Sigma')$. Because $C'^2 > 4(\gamma + 1)$
we can apply Serrano's Theorem~\ref{thm_serrano} to obtain the existence
of an effective divisor $V$ on $\text{Tor}(\Sigma')$ for which
\begin{equation} \label{serranostronger}
  0 < s < C' \cdot V - s \leq \gamma + 1 \qquad \text{and} \qquad C'^2 \leq \frac{(\gamma + 1 + s)^2}{s},
\end{equation}
where $s = V^2$. 
The bounds on $C'^2$ imply that $s = 1$, except possibly if $\gamma \in \{ 4, \dots, 13 \}$
in which case $s \in \{1, 2\}$.

We claim that this implies $h^0(\text{Tor}(\Sigma'), V) \leq s + 1$. Suppose not, and let $\Gamma$ be 
as in the proof of Theorem~\ref{maintheorem_gonality}, i.e., it is a lattice polygon 
containing at least $s + 2$ lattice points, with the property 
that $C' \cdot V \geq \text{MV}(\Delta^\text{max},\Gamma)$. 
\begin{itemize}
\item
If $\Gamma$ contains
a line segment of integral length $2$, then as before it follows 
that $C' \cdot V \geq 2\gamma$, which contradicts $C' \cdot V \leq \gamma + 1 + s$ (note that $s=1$ in case $\gamma = 3$).
\item So we can assume that $\Gamma$ does not contain such a line segment.
\begin{itemize}
\item
If $s = 1$ it therefore suffices to consider the case where $\Gamma$ contains $\Sigma$ (after performing a unimodular transformation if needed). We again
find $C' \cdot V \geq d$ where $d \geq 0$ is the smallest integer such that $\Delta^\text{max}$ is contained in a translate of $d\Sigma$.
 By definition of the lattice size, it follows that
 \[ \qquad \qquad \quad \quad C' \cdot V \geq \lsize(\Delta^\text{max}) \geq \lsize(\Delta^{(1)}) + 3 \geq \text{lw}(\Delta^{(1)}) + 5 \geq \text{lw}(\Delta) + 3 \geq \gamma + 3.  \]
Here the second inequality follows from \cite[Eq.\,(2)]{latticesize}, the third inequality follows from (\ref{planecurvecondition}),
and the fourth inequality follows from Lemma~\ref{latticewidthproperties}.(i).
This contradicts that $C' \cdot V \leq \gamma + 1 + s = \gamma + 2$.
\item
If $s = 2$ then we can assume that $\gamma \geq 4$ and that $\Gamma$ contains a unimodular 
copy of either $\square$ or $\Upsilon$. As before we respectively find that 
$C' \cdot V \geq 2\gamma$ and $C' \cdot V \geq \frac{3}{2}\gamma$. In the former case
this contradicts $C' \cdot V \leq \gamma + 1 + s$. In the latter case, the contradiction follows for $\gamma \geq 7$ only. 
To deal with the case where $\gamma \leq 6$, note that
(\ref{willbesharpened}) can be rewritten as
\[
    \qquad \qquad \quad \quad 2d = \text{lw}(d \Upsilon) > \text{lw}(\Delta^\text{max}) = \gamma, 
\]
where the last equality follows from Corollary~\ref{cacoconjecture}, and the strict inequality
in the middle holds because the lattice width of a strict subpolygon of $d\Upsilon$ is strictly less than $2d$
(we excluded the possibility that $\Delta^\text{max} \cong 2\Upsilon, 3 \Upsilon$ in the \'enonc\'e).
It follows that for $\gamma \leq 6$, the bound $C' \cdot V \geq \frac{3}{2}\gamma$ can be refined to 
$C' \cdot V \geq \frac{3}{2}(\gamma + 1)$, which is now sufficient to contradict $C' \cdot V \leq \gamma + 1 + s$.
\end{itemize}
\end{itemize}
So we conclude that indeed
$h^0(\text{Tor}(\Sigma'), V) \leq s + 1$. 
As in the proof of Theorem~\ref{maintheorem_gonality}, along with $s \leq 2$ 
this again implies that $h^0(\text{Tor}(\Sigma'), V + K) = 0$.
The remainder of the proof is an exact copy of
the corresponding part of the proof of Theorem~\ref{maintheorem_gonality} (except in the last
paragraph, where now $F \cdot C' = \gamma + 1$, but this doesn't affect the argument).
Remark that for this part we need
$g^1_{\gamma + 1}$ to be complete, which is true because the contrary would lead to infinitely
many $g^1_\gamma$'s, contradicting Corollary~\ref{numberofpencils}.

It remains to deal with the case where $\Delta^\text{max}$ is among $\Delta_1, \Delta_2, \Delta_3$. 
Here (\ref{serranostronger}) only allows us to conclude $s \in \{1,2,3\}$.
If $s \in \{1,2\}$ then the above proof applies, so we
can assume $s = 3$. We claim that in this case
$h^0(\text{Tor}(\Sigma'),V) \leq 3$.  Suppose not, then
there exists a lattice polygon $\Gamma$ containing at least $4$ lattice points,
with the property that $C' \cdot V \geq \text{MV}(\Delta^\text{max},\Gamma)$. 
\begin{itemize}
 \item If $\Gamma$ contains a line segment of integral length $2$ then we again run into
 a contradiction (note that we only consider $\gamma = 5$ and $\gamma = 6$).
 \item If not then we can again assume that $\square \subset \Gamma$ or $\Upsilon \subset \Gamma$. In the former
 case the bound $C' \cdot V \geq 2\gamma$ suffices to run into contradiction (again using that $\gamma = 5,6$).
 In the case $\Upsilon \subset \Gamma$, the above sharpened 
 bound $C' \cdot V \geq \frac{3}{2}(\gamma + 1)$ results in a contradiction for $\Delta_2$ and
 $\Delta_3$, but remains insufficient in the case of $\Delta_1$.
 Now it is not hard to see 
 that there is no unimodular transformation mapping $\Delta_1$ inside $3\Upsilon$. Indeed, because
 the lattice width of a subpolygon of $3\Upsilon$ that misses two vertices of $3\Upsilon$ is at most $4$, we find
 that a unimodular copy of $\Delta_1$ inside $3\Upsilon$ should have an edge in common with $3\Upsilon$.
 But $\Delta_1$ contains only one edge having $4$ lattice points, and the width of
 $\Delta_1$ with respect to the direction of this edge is $8$. 
  So $\Delta_1$ can indeed impossibly fit inside $3\Upsilon$.
 It follows that the smallest multiple of $\Upsilon$ containing a unimodular copy of
 $\Delta_1$ is $4\Upsilon$, from which \[ C' \cdot V  \geq \text{MV}(\Delta^\text{max},\Gamma)
  \geq \text{MV}(\Delta^\text{max},\Upsilon) \geq 3 \cdot 4 = 12. \] This gives the desired contradiction.
\end{itemize}
So we conclude that indeed
$h^0(\text{Tor}(\Sigma'), V) \leq 3$.
This implies that $h^0(\text{Tor}(\Sigma'),V + K) = 0$, and the rest of the argument
can again be copied word by word, essentially.
\hfill $\blacksquare$\\

\noindent \emph{Remark.} Kawaguchi's proof technique should in principle
allow one to obtain similar theorems on base-point free $\gamma^1_{\gamma + n}$'s for $n=2, 3, \dots$\
Here condition (\ref{planecurvecondition}) will have to be replaced by
 \[
  \lsize(\Delta^{(1)}) \geq \text{lw}(\Delta^{(1)}) + n + 1.
 \]
However, an increasing number of exceptional polygons are expected to come into play, both for geometric
reasons (definitely, more and more multiples of $\Upsilon$ will show up) and for
proof-technical reasons (as in the case of $\Delta_1, \Delta_2, \Delta_3$ in the above proof).
This might be feasible for $n=2$, although we did not try this in detail. For higher values of $n$ we expect a complete classification 
to become very complicated.




\section{Clifford index and Clifford dimension} \label{section_clifford}

To a smooth projective curve
$C / k$ of genus $g \geq 4$ one can associate its Clifford index
\[ \text{ci}(C) = \min \{ \, d - 2r \, | \, \text{$C$ carries a divisor $D$ with $|D| = g^r_d$} \qquad \qquad \qquad \qquad \qquad \]
 \[ \qquad \qquad \qquad \qquad \qquad \qquad \qquad \qquad \qquad \text{and $h^0(C,D),h^0(C,K-D) \geq 2$} \, \} \]
(where $K$ is a canonical divisor on $C$)
and its Clifford dimension
\[ \text{cd}(C) = \min \{ \, r \, | \, \text{there exists a $g^r_d$ realizing $\text{ci}(C)$} \}; \]
see \cite{ELMS}.
In the case of a singular and/or non-complete curve $C/k$, we define $\text{ci}(C)$ and $\text{cd}(C)$
to be the corresponding quantities associated to its smooth complete model.
In this section we give a combinatorial interpretation for the Clifford index and the
Clifford dimension.
Again the key trick is due to Kawaguchi \cite[Proof of Thm.\,1.3.(iii)]{Kawaguchi},
but thanks to our more careful analysis of the planar curve case we obtain a complete statement.
\begin{theorem} \label{cliffordtheorem}
   Let $f \in k[x^{\pm 1},y^{\pm 1}]$ be non-degenerate with respect to its Newton polygon
 $\Delta = \Delta(f)$ and suppose that $\sharp(\Delta^{(1)} \cap \mathbb{Z}^2) \geq 4$. Then
 \begin{itemize}
   \item if $\Delta^{(1)} \cong (d-3)\Sigma$ for $d \geq 5$ then $\emph{ci}(U(f)) = d-4$ and $\emph{cd}(U(f)) = 2$,
   \item if $\Delta^{(1)} \cong \Upsilon$ then $\emph{ci}(U(f)) = 1$ and $\emph{cd}(U(f)) = 1$,
   \item if $\Delta^{(1)} \cong 2\Upsilon$ then $\emph{ci}(U(f)) = 3$ and $\emph{cd}(U(f)) = 3$,
   \item in all other cases $\emph{ci}(U(f)) = \emph{lw}(\Delta^{(1)})$ and $\emph{cd}(U(f)) = 1$.
 \end{itemize}
\end{theorem}

\noindent \textsc{Proof.}
The first three cases correspond to smooth projective plane curves of degree $d \geq 5$, non-hyperelliptic
curves of genus $4$, resp.\
smooth intersections of pairs of cubics in $\mathbb{P}^3$,
while the cases $\Delta^{(1)} \cong \Gamma^5_1, \Gamma^5_2, \Gamma^5_3$ correspond
to non-hyperelliptic, non-trigonal curves of genus $5$.
In these 
situations the Clifford index and the Clifford dimension are well-known; see \cite[pp.\,174-175]{ELMS}
and \cite[p.\,225]{eisenbud}.
In all other cases Corollary~\ref{numberofpencils} yields that the number of gonality pencils is finite,
while from Corollary~\ref{cacoconjecture} we know that $\gamma(U(f))=
\text{lw}(\Delta^{(1)}) + 2$.
A result by Coppens and Martens \cite{CoppensMartens} (see the discussion preceding 
\cite[Thm.\,B]{CoppensMartens}) then implies that $\text{ci}(U(f)) = \text{lw}(\Delta^{(1)})$.
By definition of the Clifford dimension, this implies $\text{cd}(U(f)) = 1$. 
\hfill $\blacksquare$\\

\noindent \emph{Remark.} For curves $C/k$ of genus $1 \leq g \leq 3$ one sometimes defines
\begin{itemize}
 \item $\text{ci}(C) = 1$ if $C$ is a non-hyperelliptic genus $3$ curve, and $\text{ci}(C) = 0$ if not,
 \item $\text{cd}(C) = 1$.
\end{itemize}
With these conventions, Theorem~\ref{cliffordtheorem} remains valid when one replaces
the condition $\sharp(\Delta^{(1)} \cap \mathbb{Z}^2) \geq 4$ with $\sharp(\Delta^{(1)} \cap \mathbb{Z}^2) \geq 1$.

\begin{corollary} \label{smoothplanecorollary}
   Let $f \in k[x^{\pm 1},y^{\pm 1}]$ be non-degenerate with respect to its (two-dimen\-sional) Newton polygon
 $\Delta = \Delta(f)$. Then $U(f)$ is birationally equivalent to a smooth projective plane curve if and only
 if $\Delta^{(1)} = \emptyset$ or $\Delta^{(1)} \cong (d-3)\Sigma$ for some integer $d \geq 3$.
\end{corollary}

\noindent \textsc{Proof.} The `if' part is easily verified. As for the `only if' part, let $g$
be the geometric genus of $U(f)$, which is necessarily of the form
$(d-1)(d-2)/2$ for some $d \geq 2$. If $d \geq 5$ then $\text{cd}(U(f)) = 2$ and
the corollary follows from Theorem~\ref{cliffordtheorem}.
If $d=2$ or $d=3$ then the statement is trivial. If $d=4$ then the claim follows
because $U(f)$ is non-hyperelliptic, and because $\Sigma$ is the only
two-dimensional lattice polygon containing $g=3$ lattice points (up to unimodular equivalence). \hfill $\blacksquare$

\section{Scrollar invariants} \label{section_maroni}

We begin by recalling some facts on
rational normal scrolls and on scrollar invariants. Our main references
are \cite{EiHa}, \cite[\S8.26-29]{harris} and \cite[\S1-4]{Schreyer}.

Let $n \in \mathbb{Z}_{\geq 1}$ and let $\mathcal{E}=\mathcal{O}(e_1)\oplus \cdots \oplus \mathcal{O}(e_n)$ be a locally free sheaf of rank $n$ on $\mathbb{P}^1$. Denote by $\pi:\mathbb{P}(\mathcal{E})\to\mathbb{P}^1$ the corresponding $\mathbb{P}^{n-1}$-bundle. We assume that $0\leq e_1\leq e_2\leq \ldots\leq e_n$ and that $e_1 + e_2 + \dots + e_n \geq 2$. Set $N=e_1+e_2+\ldots+e_n+n-1$. 
A \emph{rational normal scroll of type $(e_1, \dots, e_n)$} in $\mathbb{P}^N$ is the image of
the induced morphism
$$\mu:\mathbb{P}(\mathcal{E})\to \mathbb{P}H^0(\mathbb{P}(\mathcal{E}),\mathcal{O}_{\mathbb{P}(\mathcal{E})}(1)),$$
composed with an isomorphism $\mathbb{P}H^0(\mathbb{P}(\mathcal{E}),\mathcal{O}_{\mathbb{P}(\mathcal{E})}(1)) \rightarrow \mathbb{P}^N$. 

The dimension of a rational normal scroll of type $(e_1, \dots, e_n)$ equals $n$, while its degree equals $e_1+\ldots+e_n=N-n+1$.
This means that
the classical lower bound $\deg(X)\geq \text{codim}_{\mathbb{P}^N}(X)+1$ for projective varieties $X\subset \mathbb{P}^N$ that are not contained in any hyperplane is attained. Varieties for which this holds are said to \emph{have minimal degree}.
They have been classified by Del Pezzo (the surface case, 1886) and Bertini (1907):
any projective variety of minimal degree is a cone over a smooth such variety, and the smooth such varieties are exactly the rational normal scrolls with $e_1 > 0$, the quadratic hypersurfaces, and the Veronese surface in $\mathbb{P}^5$. See \cite{EiHa} for a modern proof.

There is an easy geometric way of describing rational normal scrolls. Consider linear subspaces $\mathbb{P}^{e_1},\ldots,\mathbb{P}^{e_n}\subset\mathbb{P}^N$ that span $\mathbb{P}^N$. 
In each $\mathbb{P}^{e_\ell}$, take a rational normal curve\footnote{If $e_\ell = 0$ then this
`curve' is just a point, in fact. We will keep making this abuse of language.} of degree $e_\ell$, e.g.\ parameterized by
\begin{equation} \label{veronese}
 \nu_\ell : \mathbb{P}^1 \rightarrow \mathbb{P}^{e_\ell}: (X : Z) \mapsto \left(Z^{e_\ell} : XZ^{e_\ell - 1} : \dots : X^{e_\ell} \right) .
\end{equation}
Then $$S = \bigcup_{P\in\mathbb{P}^1} \langle \nu_1(P),\ldots,\nu_n(P)\rangle \subset \mathbb{P}^N$$ is a rational normal scroll of type $(e_1, \dots, e_n)$, and conversely every rational normal scroll arises in this way.
The scroll is smooth if and only if $e_1>0$. In this case $\mu : \mathbb{P}(\mathcal{E}) \rightarrow S$ is an isomorphism.
If $0=e_1=\ldots=e_\ell<e_{\ell+1}$ with $1 \leq \ell<n$, then the scroll is a cone with an $(\ell-1)$-dimensional vertex.
In this case $\mu : \mathbb{P}(\mathcal{E}) \rightarrow S$ is a resolution of singularities.
Outside the exceptional locus, our $\mathbb{P}^{n-1}$-bundle $\pi : \mathbb{P}(\mathcal{E}) \rightarrow \mathbb{P}^1$ corresponds to
\[
 S \setminus S^\text{sing} \to\mathbb{P}^1: Q\in\langle \nu_1(P),\ldots,\nu_n(P)\rangle \mapsto P.
\]
Abusing notation, we denote this map also by $\pi$. Abusing terminology, when talking about the fiber of $\pi$ above a point $P$, we
mean the whole space $\langle \nu_1(P),\ldots,\nu_n(P)\rangle$.

Now let $C\subset \mathbb{P}^{g-1}$ be a canonical curve of genus $g \geq 3$ 
and fix any pencil $g_d^1$ on $C$. Let $K \sim \mathcal{O}_C(1)$ be a 
canonical divisor on $C$. For an effective divisor
$D \in g^1_d$, denote by $\langle D \rangle$ its linear span 
(if $D$ is the sum of $\gamma$ distinct points, the linear span of $D$ is just the
linear span of these points; in general one defines it as the intersection of all hyperplanes
whose intersection divisor with $C$ is at least $D$, see \cite[\S2.3]{Schreyer}).
The Riemann-Roch theorem implies that $h^0(C, K - D) = g - d - 1 + h^0(C,D)$, from which it follows that
the dimension of $\langle D \rangle$ equals
\begin{equation} \label{completeornot}
d - h^0(C,D).
\end{equation}
This does not depend the specific choice of $D$. In particular, if our $g^1_d$ is complete, then the dimension of $\langle D \rangle$ is $d - 2$.

Consider
\begin{equation} \label{formulaforscroll}
 S=\bigcup_{D\in g_d^1}\, \langle D\rangle \subset \mathbb{P}^{g-1}.
\end{equation}
Then $S$ is a rational normal scroll by \cite[Thm.\,2]{EiHa} or \cite[(2.5)]{Schreyer}, and it contains the curve $C$.
In most interesting cases $\dim S = d - h^0(C,D) +1$, but it may happen
that $\dim S = d - h^0(C,D)$, which holds iff $h^0(C, K-D) = 0$, i.e.\ iff $\langle D \rangle = \mathbb{P}^{g-1}$.
 If $g^1_d$ is base-point free then $C$ does not meet the singular locus of $S$ (in which case the restriction of $\pi$ to $C$ is
a dominant rational map of degree $d$).

Let
$(e_1, \dots, e_n)$ be the type of $S$. Then the numbers $e_1, \dots, e_n$
are called the \emph{scrollar invariants of $C$ with respect to $g^1_d$}. 
When we talk about the \emph{scrollar invariants of $C$}, without making
reference to a specific pencil, we always mean the scrollar invariants
with respect to a gonality pencil, but note that this may depend on the choice
of the latter, in which case the terminology is avoided. In the
trigonal case the notion is well-behaved, and here the scrollar invariants
are better known under the name \emph{Maroni invariants}.\footnote{
Unfortunately, the existing literature
is ambiguous at this point: sometimes one talks about \emph{the} Maroni invariant of a trigonal curve, 
in which case one could mean either $e_1$ or $e_2 - e_1$.}
The scrollar invariants of an arbitrary
non-hyperelliptic curve $C/k$ of genus $g \geq 3$ with respect to a
pencil $g^1_d$ are then defined to be the corresponding invariants of a canonical model.

If $g^1_d = |D|$ is complete and base-point free then $n = d-1$, and
the scrollar invariants can alternatively be described as follows:
$$h^0(C,mD)=\begin{cases} h^0(C,(m-1)D)+1=m+1 & \text{if } 0\leq m\leq e_1+1, \\ h^0(C,(m-1)D)+2 & \text{if } e_1+1<m\leq e_2+1, \\ \hspace{1.5cm} \vdots & \hspace{1.5cm} \vdots \\ h^0(C,(m-1)D)+ d - 1 & \text{if } e_{d-2}+1<m\leq e_{d-1}+1, \\ h^0(C,(m-1)D)+d=md-g+1 & \text{if } m>e_{d-1}+1.\end{cases}$$
See \cite[(2.4)]{Schreyer} for more details, as well as a treatment of the general case (where our $g^1_d$
is not necessarily complete and/or base-point free).\\

\noindent \emph{Remark.}
  From this description it follows that if our $g^1_d$ is complete and base-point free
  then $e_{d-1}\leq \frac{2g-2}{d}$. 
  Indeed, if $m>\frac{2g-2}{d}$ then $h^0(C,mD)=md-g+1$ 
  and
  by the above characterization, the smallest $m$ for which $h^0(C, mD)=md-g+1$ is $m=e_{d-1}+1$.\\

The main result of this section is as follows.

\begin{theorem} \label{thm_maroni}
Let $f \in k[x^{\pm 1},y^{\pm 1}]$ be non-degenerate with respect to its Newton polygon
$\Delta = \Delta(f)$, and assume that $\Delta^{(1)}$ is two-dimensional.
Let $v$ be a lattice direction. Then the multiset of scrollar invariants of $U(f)$ with respect to $g_v$ equals
the multiset of non-negative width invariants of $\Delta$ with respect to $v$.
\end{theorem}

\noindent \emph{Remark.} As mentioned at the end of Section~\ref{section_nondegistoric}, our main results
stay true if one weakens the assumption of being $\Delta$-non-degenerate to being $\Delta$-toric. 
This also applies to Theorem~\ref{thm_maroni}, but the argument becomes more technical due to the potential presence of base points. For
the sake of clarity, the proof below only handles the case of $\Delta$-non-degenerate Laurent polynomials.
The extra ingredients in the $\Delta$-toric case are then sketched in a following remark.\\

\noindent \textsc{Proof.}
Write $d = w(\Delta,v)$, so that $g_v$ is a base-point free $g^1_d$.
Using a unimodular transformation if needed, we may assume that $v = (a,b) = (1,0)$ and that
$\Delta$ is contained in the horizontal strip $\mathbb{R} \times [0,d] \subset \mathbb{R}^2$.
Then the width invariants of $\Delta$ with respect to $v$ are the numbers
\[ E_\ell=\sharp\{(i,j)\in \Delta^{(1)}\cap\mathbb{Z}^2\,|\,j=\ell\} \, - \, 1, \]
where $\ell = 1, \dots, d - 1$. We have to show that the scrollar
invariants with respect to the pencil cut out by $p : U(f) \mapsto \mathbb{T}^1 : (x,y) \mapsto x$ are given by the multiset
$\{ E_\ell \}_{\ell = 1, \dots, d-1} \cap \mathbb{Z}_{\geq 0}$. Denote the cardinality of
this multiset by $n$.

Let $C$ be the canonical model of $U(f)$ obtained by taking the Zariski closure of its image under
the morphism $\varphi_{\Delta^{(1)}}$, as described in (\ref{canonicalembedding}).
For all $\ell\in\{1,\ldots,d -1\}$ for which $E_\ell \geq 0$, 
let $\mathbb{P}^{E_\ell}\subset \mathbb{P}^{g-1}$ 
be the linear subspace defined by $X_{i,j}=0$ for all $(i,j)\in \Delta^{(1)}\cap\mathbb{Z}^2$ for 
which $j\neq \ell$. That is, $\mathbb{P}^{E_\ell}$ is the subspace corresponding to the projective 
coordinates $(X_{i,\ell})_{(i,\ell)\in\Delta^{(1)}\cap \mathbb{Z}^2}$. Also consider the rational normal 
curves parameterized by $\nu_\ell : \mathbb{P}^1 \to \mathbb{P}^{E_\ell}$
as in (\ref{veronese}), i.e.
$$\forall x \in k^\ast : \nu_\ell(x) = (1:x:\ldots:x^{E_\ell}).$$ 
Then $\varphi_{\Delta^{(1)}}$ maps every $(x,y) \in \mathbb{T}^2$ inside
the $(n-1)$-dimensional linear subspace of $\mathbb{P}^{g-1}$ spanned by the points $\nu_\ell(x)$. Indeed, abusing 
notation, one sees that when the $\nu_\ell(x)$'s
are scaled by an appropriate power of $x$, the point $\varphi_{\Delta^{(1)}}(x,y)$
arises as the linear combination
\[ \sum_{\substack{\ell = 1 \\ E_\ell \geq 0}}^{d-1} y^\ell \nu_\ell(x).\]
Now for all but finitely many $c \in k^\ast$, the inverse image divisor $p^{-1}(c)$ consists of $d$ distinct 
points $(c,y_1), \dots, (c,y_d)$ of $U(f)$.
For these $c$, the linear span $\langle D_c \rangle$ of $D_c=\varphi_{\Delta^{(1)}}(p^{-1}(c))$ is contained in 
$\langle\nu_\ell(c)\rangle_\ell$, and since the matrix
\[ \left( y_i^\ell \right)_{\substack{i=1, \dots, d \\ \ell = 1, \dots, d - 1 \\ E_\ell \geq 0}} \]
has rank $n$ (indeed, its columns are linearly independent because
by adding a number of columns one obtains a $(d \times d)$-Vandermonde matrix), we find
that actually $\langle D_c \rangle = \langle\nu_\ell(c)\rangle_\ell$.
We conclude that the scroll $S \subset \mathbb{P}^{g-1}$
swept out by our $g^1_d$ is exactly the rational normal scroll parameterized by the $\nu_\ell$'s.
Hence we obtain that the multiset of scrollar invariants with respect to $g_d^1$ equals the multiset consisting
of the non-negative $E_\ell$'s, which is exactly what we wanted. \hfill $\blacksquare$\\

\noindent \emph{Remark (continued).} If $f$ is only $\Delta$-toric, rather than $\Delta$-non-degenerate, it may happen that
$d' < d$, where $d' = w(\Delta(f),v)$ and $d = w(\Delta,v)$. 
In this case $g_v$ decomposes into a base-point free $g^1_{d'}$ and a fixed part $F$ which is supported
on at most two zero-dimensional toric orbits, as explained
in the proof of Lemma~\ref{gv_basepointfree}.
The base-point free part corresponds to the morphism $p : U(f) \rightarrow \mathbb{T}^1 : (x,y) \mapsto x$,
 and the above reasoning shows that for all but finitely many $c \in k^\ast$, the linear span $\langle D_c \rangle$
  of $D_c = \varphi_{\Delta^{(1)}}(p^{-1}(c))$ equals $\langle\nu_\ell(c)\rangle_\ell$. 
For each of these $D_c$ one clearly has $\langle D_c \rangle \subset \langle D_c + F \rangle$. We claim that actually
equality holds. This implies that the scroll swept out by $g_v$ coincides with the scroll swept out by its base-point free part, so that
Theorem~\ref{thm_maroni} also follows in the $\Delta$-toric case.
Note that it suffices to prove the claim under the assumption that $\Delta = \Delta^\text{max}$. Indeed,
from Lemma~\ref{cruciallemma} (and the consequent remark) we see that if $f$ is $\Delta$-toric, then it is also $\Delta^\text{max}$-toric.
Of course switching from $\Delta$ to $\Delta^\text{max}$ may have an influence on $g_v$, but it 
can only affect the fixed part $F$, and if it does
then $F$ becomes replaced by $F'$ with $F' > F$. So if we can prove that $\langle D_c \rangle = \langle D_c + F' \rangle$ then necessarily
$\langle D_c \rangle = \langle D_c + F \rangle$.

Let $\Delta(f)$ be 
contained in the strip $m_f \leq Y \leq m_f + d'$ and suppose
that $0 < m_f$. Recall that $\Delta$ has a unique lower-most vertex $v_\text{low}$. Our assumption
$\Delta = \Delta^\text{max} = \Delta^{(1)(-1)}$ ensures that also $\Delta^{(1)}$ has a unique lower-most vertex 
and that the adjacent cones are similar. Denote the corresponding zero-dimensional orbit by $P$.
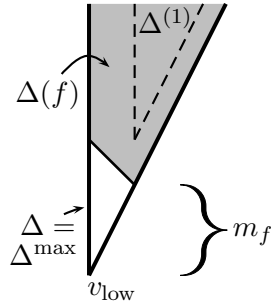
\begin{center}
\psset{unit=0.3cm}
\begin{pspicture}(0,-1)(6,12.5)
\pspolygon[linestyle=none,fillstyle=solid,fillcolor=lightgray](2,4)(6,12)(0,12)(0,6)
\psline[linewidth=1.5pt](0,0)(6,12)
\psline[linewidth=1.5pt](0,0)(0,12)
\psline[linewidth=1pt](0,6)(2,4)
\psline[linestyle=dashed](2,6)(5,12)
\psline[linestyle=dashed](2,6)(2,12)
\psbrace(4,0)(4,4){ }
\rput(7.2,1.83){\small $m_f$}
\rput(-1.8,8){\small $\Delta(f)$}
\pscurve{->}(-1.2,8.9)(-0.1,9.7)(0.9,9.2)
\rput(-1.8,2.2){\small $\Delta =$}
\rput(-2,0.9){\small $\Delta^\text{max}$}
\pscurve{->}(-1.1,2.7)(-0.6,2.9)(-0.2,3)
\rput(3.35,11.3){\small $\Delta^{(1)}$}
\rput(1,-0.7){\small $v_\text{low}$}
\end{pspicture}
\end{center}
Then locally around $O(v_\text{low})$ we have a natural isomorphism $\text{Tor}(\Delta) \rightarrow \text{Tor}(\Delta^{(1)})$
under which $O(v_\text{low})$ corresponds to $P$. From the proof of Lemma~\ref{gv_basepointfree}
we conclude that $C \subset \text{Tor}(\Delta^{(1)})$ intersects the zero divisor $\mathcal{F}_c$ of $x^ay^b - c$, with $c \in k^\ast$ sufficiently generic,
with multiplicity $m_f$ in $P$. Our task is to prove that every hyperplane $H$ containing the support of $D_c$ intersects 
$C$ in $P$ with multiplicity at least $m_f$. But this follows from
$$I_{P,\mathbb{P}^{g-1}}(H,C) = I_{P,\text{Tor}(\Delta^{(1)})}(H\cap \text{Tor}(\Delta^{(1)}),C) \geq I_{P,\text{Tor}(\Delta^{(1)})}(\mathcal{F}_c,C)=m_f,$$
where $I_{P,X}(\cdot, \cdot)$ denotes the intersection multiplicity of the arguments in $P$ when viewed as schemes inside $X$, and the inequality holds because $H \supset \langle D_c \rangle = \langle \nu_\ell(c) \rangle_\ell \supset \mathcal{F}_c$. A similar reasoning at the top (if needed) then proves
that $\langle D_c \rangle = \langle D_c + F \rangle$.\\

\noindent \emph{Example (revisited, see Sections~\ref{section_combpencils} and~\ref{section_gonality}).} Let $\Delta$ be the lattice polygon
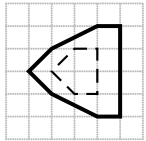
\begin{center}
\psset{unit=0.3cm}
\begin{pspicture}(-1,-3)(5,3)
\psgrid[subgriddiv=1,griddots=10,gridcolor=lightgray,gridlabels=0pt](-1,-3)(5,3)
\pspolygon[linewidth=1.5pt](0,0)(1,-1)(3,-2)(4,-2)(4,2)(3,2)(1,1)
\pspolygon[linestyle=dashed](1,0)(2,-1)(3,-1)(3,1)(2,1)
\end{pspicture}
\end{center}
where $\Delta^{(1)}$ is marked in dashed lines. Let $f \in k[x^{\pm 1}, y^{\pm 1}]$ be a $\Delta$-non-degenerate (or $\Delta$-toric) 
Laurent polynomial. Then $U(f)$ is a $4$-gonal genus $7$ curve carrying exactly
two $g^1_4$'s, namely $g_{(1,0)}$ and $g_{(0,1)}$. In the former case the scrollar invariants are $\{1,1,2\}$ while in the latter case
they read $\{ 0,2,2\}$.\\

As a corollary to the proof of Theorem~\ref{thm_maroni} we find:

\begin{corollary} \label{completerank}
 Let $f \in k[x^{\pm 1},y^{\pm 1}]$ be non-degenerate with respect to its Newton polygon
$\Delta = \Delta(f)$, where we assume that $\Delta^{(1)}$ is two-dimensional.
Let $v$ be a lattice direction. Then the rank of the complete linear system
spanned by $g_v$ equals the number of negative width invariants 
of $\Delta$ with respect to $v$ (counting multiplicities) plus $1$.
\end{corollary}

\noindent \textsc{Proof.} Let $d = w(\Delta,v)$ 
and let $D \in g_v$, and assume that we work
on the canonical model $C$ of $U(f)$ from the proof of Theorem~\ref{thm_maroni}.
By (\ref{completeornot}) we
know that $\langle D \rangle$ is $(d - h^0(C,D))$-dimensional,
while the proof of Theorem~\ref{thm_maroni} tells us that the dimension equals
the number of non-negative lattice width invariants minus $1$. From this
the statement follows. \hfill $\blacksquare$\\

In particular we find the following combinatorial characterization of completeness:

\begin{corollary} \label{nonnegativewidth}
 Let $f \in k[x^{\pm 1},y^{\pm 1}]$ be non-degenerate with respect to its Newton polygon
$\Delta = \Delta(f)$, where we assume that $\Delta^{(1)}$ is two-dimensional.
Let $v$ be a lattice direction. Then $g_v$ is complete if and only if
the width invariants of $\Delta$ with respect to $v$ are all non-negative.
\end{corollary}

\noindent \emph{Example.} Let $\Delta = d\Sigma$ for some $d \geq 2$, so
that $U(f)$ is birationally equivalent to a smooth plane curve of degree $d$.
The width invariants of $\Delta$ with respect to $(1,0)$ are $(d-3, d-4, \dots, 1, 0, -1)$, so
that $g_{(1,0)}$ is not complete. (Indeed: it is a subsystem of the $g^2_d$ cut out
by all line sections of $\mathbb{P}^2$.)\\

\begin{corollary} \label{nonnegativewidth}
 Let $f \in k[x^{\pm 1},y^{\pm 1}]$ be non-degenerate with respect to its Newton polygon
$\Delta = \Delta(f)$, where we assume that $\Delta^{(1)}$ is two-dimensional.
Then the dimension of the scroll spanned by $g_v$ equals the number of non-negative lattice width invariants
of $\Delta$ with respect to $v$,
unless this number is $g$ (i.e.\ there are no strictly positive lattice width invariants) in which case the dimension equals $g - 1$.
\end{corollary}

\noindent \textsc{Proof.} This follows from the considerations below formula (\ref{formulaforscroll}), along with
the combinatorial interpretation for $d - h^0(C,D)$ stated in the proof of Corollary~\ref{completerank}. \hfill $\blacksquare$\\

\vspace{0.2cm}
\noindent \emph{Remarks.}
\begin{itemize}
\item
Inheriting the notation of the proof of Theorem~\ref{thm_maroni}, we have $C \subset \text{Tor}(\Delta^{(1)}) \subset S \subset \mathbb{P}^{g-1}$.
One can verify that $\text{Tor}(\Delta^{(1)})$ intersects the fiber
      of $\pi$ above a point $x \in k^\ast$ in a rational normal curve of degree $\gamma-2$.
      Above $(1:0),(0:1) \in \mathbb{P}^1$ this fiber may degenerate.
\item Through Corollary~\ref{cacoconjecture} and Theorem~\ref{thm_maroni}, 
the upper bound $\frac{2g-2}{\gamma}$ on the scrollar invariants
with respect to a gonality pencil $g^1_\gamma$
implies the purely combinatorial inequality
\[ \text{lw}(\Delta) \cdot E_\ell \leq 2\,  \sharp (\Delta^{(1)} \cap \mathbb{Z}^2) - 2, \]
where $\Delta$ is understood to be contained in
\[ \{(i,j)\in \mathbb{R}^2\,|\, 0\leq j\leq \text{lw}(\Delta)\} \]
and the $E_\ell$'s are the width invariants of $\Delta$ with respect to any lattice width direction. This inequality
holds as soon as $\sharp (\Delta^{(1)} \cap \mathbb{Z}^2) \geq 1$ (including the cases $\Delta = 2\Upsilon$
and $\Delta = d \Sigma$, which can be verified separately).
The bound can be attained. For example, consider the lattice polygon $\Delta_{a,b}=\text{conv} \{ (b,0), (0,a), (0,0) \}$, where $a \geq 2$ and $b$ is of the form $ak-1$ with $k\in\mathbb{Z}_{\geq 2}$. In this case, $\gamma=\text{lw}(\Delta_{a,b})=a$ 
is computed by $(1,0)$, and $E_1=ak-k-2=\frac{2g-2}{\gamma}$.
\end{itemize}

\section{Applications} \label{section_applications}

\subsection*{Curves with prescribed invariants}
The results of this article
might serve as a tool in proving 
certain existence results in Brill-Noether theory.
The number of inequivalent lattice polygons grows very
quickly with the genus (for instance, in genus $30$ this number is $957\,001$; see \cite[Tab.\,1]{movingout}), 
resulting in a wide variety of Brill-Noether types, that (at least in principle) can 
be scanned by exhaustive search. To highlight one example, let $\Delta$ be the following polygon.
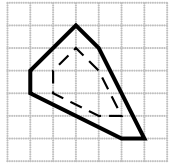
\begin{center}
\psset{unit=0.3cm}
\begin{pspicture}(-1,-1)(6,6)
\psgrid[subgriddiv=1,griddots=10,gridcolor=lightgray,gridlabels=0pt](-1,-1)(6,6)
\pspolygon[linestyle=dashed](3,1)(4,1)(3,3)(2,4)(1,3)(1,2)
\pspolygon[linewidth=1.5pt](4,0)(5,0)(3,4)(2,5)(0,3)(0,2)
\end{pspicture}
\end{center}
Every $\Delta$-non-degenerate (or $\Delta$-toric) curve is a $5$-gonal curve of
genus $9$ admitting exactly three $g^1_5$'s (corresponding to the lattice directions
$(1,0)$, $(0,1)$ and $(1,-1)$), that are independent of each other, and with respect 
to each of which the scrollar invariants are $\{ 0, 1, 2, 2\}$.
Moreover, by \cite[Thm.\,2.5.12]{Koelman} the locus of such curves 
inside the moduli space $\mathcal{M}_9$ of curves of genus $9$ has dimension $15$. 
See \cite{coppensfivegonal} for a related discussion; note
that each of our
$g^1_5$'s is of `type II' (i.e.\ $0$ is among the scrollar invariants), as opposed
to the `type I' pencils that are the main object of study in \cite{coppensfivegonal}.

We want to stress that many Brill-Noether
types are not represented in the toric world. For instance, Lemma~\ref{latticewidthproperties}.(vi) shows
that the gonality of a smooth curve in a toric surface is $O(\sqrt{g})$, while
general curves of genus $g$ have gonality $\lceil g/2 \rceil + 1$.
So the class of curves that
we are considering in this article is rather special.
In terms of moduli, the locus of curves of genus $g \geq 4$ that admit a smooth
embedding in a toric surface has dimension $2g + 1$, with the exception of $g=7$, where the dimension reads
$16$; see \cite{CastryckVoight}. Recall that $\dim \mathcal{M}_g = 3g - 3$.

\subsection*{Weierstrass semi-groups of embedding dimension $2$}
The Weierstrass semi-group of a point $P$ on a smooth projective curve $C$ is the
set of possible pole orders at $P$ of functions that are regular on $C \setminus \{P\}$.
This is a numerical semi-group, i.e.\ a sub-semi-group of $\mathbb{N}$ with finite complement. 
A numerical semi-group is said to have \emph{embedding dimension $2$} if it is of the form
$a\mathbb{N} + b\mathbb{N}$ for coprime integers $a,b \geq 2$.
Using Corollary~\ref{cacoconjecture} we can prove the following:

\begin{theorem} \label{Cabtheorem}
If a smooth projective curve $C/k$ carries a point $P$ having a
Weierstrass semi-group of embedding dimension $2$, then this semi-group
does not depend on the choice of $P$.
\end{theorem}

\noindent \emph{Remark.} This is well-known in the case of hyperelliptic curves
of genus $g \geq 2$,
all of whose Weierstrass points have semi-group $2 \mathbb{N} + (2g+1)\mathbb{N}$.\\

\noindent \textsc{Proof.} If $C$ has a Weierstrass point with semi-group $a\mathbb{N} + b\mathbb{N}$
for coprime integers $a,b \geq 2$, then it is of genus $(a-1)(b-1)/2$ (by Riemann-Roch -- this is the number
of gaps in the semi-group). We claim that $C$ has gonality $\min \{a,b\}$. Together, this implies
that $a$ and $b$ are indeed uniquely determined (up to order). To prove the claim,
we use a result of Miura stating that $C$ is birationally equivalent to a smooth affine
curve of the form
\[ c_{b,0}x^b + c_{0,a}y^a + \sum_{ia + jb < ab} c_{i,j}x^iy^j, \qquad c_{b,0}c_{0,a} \neq 0. \]
See \cite[Thm.\,5.17, Lem.\,5.30]{Miura} or \cite{Matsumoto}.
From this it is clear
that $C$ is
$\Delta_{a,b}$-toric, where
\[ \Delta_{a,b} = \text{conv} \{ (b,0), (0,a), (0,0) \}\]
(in fact $C$ is even $\Delta_{a,b}$-non-degenerate, since an affine translation ensures appropriate behavior with respect to the toric boundary).
By Corollary~\ref{cacoconjecture}, we have that the gonality of $C$ equals $\text{lw}(\Delta_{a,b}) = \min\{a,b\}$. \hfill $\blacksquare$\\

\noindent \emph{Remark.} Miura studied curves having a Weierstrass semi-group of the form $a\mathbb{N} + b\mathbb{N}$
in the context of coding theory; he called them \emph{$C_{a,b}$ curves}. (In a recent past, $C_{a,b}$ curves have enjoyed fair
interest from researchers in explicit algebraic geometry \cite{DVCab,cantor,Miura}).
Then another way to state Theorem~\ref{Cabtheorem} is that a curve cannot be simultaneously $C_{a,b}$ and $C_{a',b'}$
for distinct pairs $\{a,b\}$ and $\{a',b'\}$.\\
\vspace{-0.1cm}

\subsection*{Curves in Hirzebruch surfaces} 
We can use Theorem~\ref{thm_maroni} to compute the scrollar invariants
of smooth curves on Hirzebruch surfaces.
An immediate corollary to this computation is that
if a non-hyperelliptic
   smooth projective curve $C$ of genus $g \geq 2$ can be embedded in the $n^\text{th}$ Hirzebruch surface $\mathcal{H}_n$,
   then $n$ is actually an invariant of $C$ (that is, it cannot be embedded
   in $\mathcal{H}_{n'}$ for $n' \neq n$).
\begin{theorem}
  \begin{itemize}
    \item The scrollar invariants (with respect to any gonality pencil) of a smooth projective plane curve $C/k$ of degree $d \geq 4$
    are $\{0,1, \dots, d-3\}$.
    \item The scrollar invariants (with respect to any gonality pencil) of a smooth projective curve $C/k$ of genus $g \geq 2$ and gonality $\gamma$ in
    the $n^\text{th}$ Hirzebruch surface $\mathcal{H}_n$ are
    \[ \left\{ \frac{g}{\gamma - 1} + \left(\ell - \frac{\gamma}{2} \right) n - 1 \right\}_{1 \leq \ell \leq \gamma - 1}.\]
    In particular, if $\gamma > 2$ then
    \[ n = \frac{2g - 2(\gamma - 1)(e_1 + 1)}{(\gamma - 1)(\gamma - 2)} \]
    is an invariant of the curve.
  \end{itemize}
\end{theorem}
\noindent \textsc{Proof.}
Because $\mathcal{H}_1$ is a blow-up of $\mathbb{P}^2$, the first statement is actually a corollary to the second.
Nevertheless, we will treat it separately.

Let $C \subset \mathbb{P}^2$ be a smooth projective curve of degree $d$ and fix a gonality pencil $g^1_{d-1}$ on $C$.
By \cite[Prop.\,3.13(ii)]{Serrano}, the latter is computed by projecting from a point of the curve. Using a projective transformation
we may assume that this point is $(0:1:0)$. 
Let $F(X,Y,Z)$ be a corresponding
defining homogeneous degree $d$ polynomial.
Then $F(x,y,1)$ is $\Delta$-toric, with
\[ \Delta = \text{conv} \{ (0,0), (d, 0), (1,d-1), (0,d-1) \}, \]
and our $g^1_{d-1}$ corresponds to $(x,y) \mapsto x$, i.e.\ it equals $g_{(1,0)}$. The statement now follows from Theorem~\ref{thm_maroni}.

Next, let $C$ be a smooth projective curve in $\mathcal{H}_n$.
Due to the toric description of Hirzebruch surfaces \cite[Ex.\,3.1.16]{coxlittleschenck}
we may assume that our curve $C$ is $\Delta$-toric, with
\[ \Delta = \text{conv} \{ (0,0), (a + dn, 0), (a,d), (0,d) \} \]
for integers $a \in \mathbb{Z}_{\geq 0}$ and $d \in \mathbb{Z}_{\geq 2}$.
Now
\begin{itemize}
 \item If $a=0$ and $n=1$ then $C$ is isomorphic to a smooth projective plane curve (of degree $d$)
and the statement follows from the first part.
 \item If $a > 0$ or $n > 1$ then by Theorem~\ref{maintheorem_gonality} there exists
 only one gonality pencil, corresponding to vertical projection (i.e.\ $\gamma = d$).
One finds that
\[ g = \sharp(\Delta^{(1)} \cap \mathbb{Z}^2) = \frac{\gamma(\gamma - 1)}{2} n + (\gamma - 1)(a - 1) \]
and, by Theorem~\ref{thm_maroni}, 
\[ e_\ell = a - 2 + \ell n \text{ (for $1 \leq \ell \leq \gamma - 1$)}.\]
From these two equalities the statement follows.
 \item If $n=0$ then $\Delta = [0,a] \times [0,d]$ is a standard rectangle. If $a \neq d$ then
 by Theorem~\ref{maintheorem_gonality} there exists only one gonality pencil, corresponding to
 horizontal or vertical projection (i.e.\ $\gamma = d$ or $\gamma = a$). If $a = d$ then there are two gonality pencils. In both cases
 the statement follows similarly from Theorem~\ref{thm_maroni}. \hfill $\blacksquare$
\end{itemize}

\section*{Acknowledgments}
We thank Marc Coppens, C\'edric P\'epin, Josef Schicho, Frank-Olaf Schreyer, Karl Schwede, the referee of a prior
submission, and the referees of the current submission for several helpful suggestions, corrections and/or discussions.
We also express our gratitude to Ryo Kawaguchi for sending
us a preprint of \cite{Kawaguchi}, which was the main source of inspiration for this research.
The first author thanks the Massachusetts Institute of Technology for its hospitality.
This work was partially supported by research project G093913N of FWO-Vlaanderen.

\vspace{0.5cm}
\small

\noindent \textsc{Vakgroep Wiskunde, Universiteit Gent}\\
\noindent \textsc{Krijgslaan 281, 9000 Gent, Belgium}\\
\vspace{-0.4cm}

\noindent \textsc{Departement Elektrotechniek, KU Leuven}\\
\noindent \textsc{Kasteelpark Arenberg 10/2452, 3001 Leuven, Belgium}\\
\vspace{-0.4cm}

\noindent \emph{E-mail address:} \verb"wouter.castryck@gmail.com"\\
\vspace{-0.2cm}

\noindent \textsc{Department of Mathematics and Applied Mathematics, University of Cape Town}\\
\noindent \textsc{Private Bag X1, Rondebosch 7701, South Africa}\\
\vspace{-0.4cm}

\noindent \emph{E-mail address:} \verb"filip.cools@uct.ac.za"\\


\begin{thebibliography}{99}
 \bibitem{Arbarello} Enrico Arbarello, Maurizio Cornalba, Phillip Griffiths, Joe Harris, \emph{Geometry of algebraic curves, vol.\ I}, Grundlehren der mathematischen Wissenschaften \textbf{267}, Springer (1985)
 \bibitem{baker} Matthew Baker, \emph{Specialization of linear systems from curves to graphs}, Algebra \& Number Theory \textbf{2}(6), pp.\ 613-653 (2008)
 \bibitem{batyrev} Victor Batyrev, \emph{Variations of the mixed Hodge structure of affine hypersurfaces in algebraic
tori}, Duke Mathematical Journal \textbf{69}(2), pp.\ 349-409 (1993)
\bibitem{beelen} Peter Beelen, \emph{A generalization of Baker's theorem}, Finite Fields and Their
Applications \textbf{15}(5), pp.\ 558-568 (2009)
 \bibitem{Brunsetal} Winfried Bruns, Joseph Gubeladze, Ng\^{o} Vi\^{e}t Trung, \emph{Normal polytopes, triangulations, and Koszul algebras}, Journal f\"ur die reine und angewandte Mathematik \textbf{405}, pp.\ 123-160 (1997) 
 \bibitem{movingout} Wouter Castryck, \emph{Moving out the edges of a lattice polygon}, Discrete and Computational Geometry \textbf{47}(3), pp.\ 496-518 (2012)
\bibitem{CaCoCanonical} Wouter Castryck, Filip Cools, \emph{A minimal set of generators for the canonical ideal of a non-degenerate curve}, Journal of the 
Australian Mathematical Society \textbf{98}(3), pp.\ 311-323 (2015)
 \bibitem{CaCo} Wouter Castryck, Filip Cools, \emph{Newton polygons and curve gonalities}, Journal of Algebraic Combinatorics \textbf{35}(3), pp.\ 345-366 + err.\ pp.\ 367-372 (2012)
 \bibitem{latticesize} Wouter Castryck, Filip Cools, \emph{The lattice size of a lattice polygon}, Journal of Combinatorial Theory, Series A \textbf{136}, pp.\ 64-95 (2015)
 \bibitem{CDV} Wouter Castryck, Jan Denef, Fr\'ederik Vercauteren, \emph{Computing zeta functions of nondegenerate curves},
 International Mathematics Research Papers Vol.\ \textbf{2006}, Article ID 72017, pp.\ 1-57 (2006)
 \bibitem{CastryckVoight} Wouter Castryck, John Voight, \emph{On nondegeneracy of curves}, Algebra \& Number Theory \textbf{3}(3), pp.\ 255-281 (2009)
 \bibitem{Coppens} Marc Coppens, \emph{The number of linear systems computing the gonality}, Journal of the Korean Mathematical Society \textbf{37}(3), pp.\ 437-454 (2000)
 \bibitem{coppensfivegonal} Marc Coppens, \emph{Five-gonal curves of genus nine}, Collectanea Mathematica \textbf{56}(1), pp.\ 21-26 (2005)
 \bibitem{CoppensMartens} Marc Coppens, Gerriet Martens, \emph{Secant spaces and Clifford's theorem}, Compositio Mathematica \textbf{78}(2), pp.\ 193-212 (1991)
 \bibitem{coxlittleschenck} David Cox, John Little, Hal Schenck, \emph{Toric varieties}, Graduate Studies in Mathematics \textbf{124}, American Mathematical Society (2011)
 \bibitem{DVCab} Jan Denef, Fr\'ederik Vercauteren, \emph{Computing zeta functions of $C_{ab}$ curves using Monsky-Washnitzer cohomology}, Finite Fields and Their Applications \textbf{12}(1), pp.\ 78-102 (2006)
 \bibitem{draismaetal} Jan Draisma, Tyrrell McAllister, Benjamin Nill, \emph{Lattice width directions and Minkowski's $3^d$-theorem}, SIAM Journal on Discrete Mathematics \textbf{26}(3), pp.\ 1104-1107 (2012)
 \bibitem{eisenbud} David Eisenbud, \emph{The geometry of syzygies}, Graduate Texts in Mathematics \textbf{229}, Springer (2005)
 \bibitem{EiHa} David Eisenbud, Joe Harris, \emph{On Varieties of Minimal Degree (A Centennial Account)}, Proceedings in Symposia in Pure Mathematics \textbf{46}, pp.\ 3-13 (1987)
 \bibitem{ELMS} David Eisenbud, Herbert Lange, Gerriet Martens, Frank-Olaf Schreyer, \emph{The Clifford dimension of a projective curve}, Compositio Mathematica \textbf{72}(2), pp.\ 173-204 (1989)
 \bibitem{FejesMakai} L\'aszl\'o Fejes T\'oth, Endre Makai Jr., \emph{On the thinnest non-separable lattice of convex
plates}, Studia Scientiarum Mathematicarum Hungarica \textbf{9}, pp.\ 191-193 (1974)
 \bibitem{fulton} William Fulton, \emph{Introduction to toric varieties}, Annals of Mathematics Studies \textbf{131}, Princeton University Press (1993)
 \bibitem{GKZ} Israel Gelfand, Mikhail Kapranov, Andrei Zelevinsky, \emph{Discriminants, resultants, and multidimensional determinants}, Birkh\"auser Boston (1994)
 \bibitem{HaaseSchicho} Christian Haase, Josef Schicho, \emph{Lattice polygons and the number $2i+7$}, American Mathematical Monthly \textbf{116}(2), pp.\ 151-165 (2009)
 \bibitem{cantor} Ryuichi Harasawa, Joe Suzuki, \emph{Fast Jacobian group arithmetic on $C_{ab}$ curves}, Proceedings of ANTS-IV (Leiden, The Netherlands), Lecture Notes in Computer Science \textbf{1838}, pp.\ 359-376 (2000)
 \bibitem{harris} Joe Harris, \emph{Algebraic geometry: a first course}, Graduate Texts in Mathematics \textbf{133}, Springer (1992)
 \bibitem{hartshorne} Robin Hartshorne, \emph{Algebraic geometry}, Graduate Texts in Mathematics \textbf{52}, Springer (1977)
 \bibitem{Harui} Takeshi Harui, \emph{The gonality and the Clifford index of curves on an elliptic ruled surface}, Archiv der Mathematik \textbf{84}, pp.\ 131-147 (2005)
 \bibitem{Kawaguchi} Ryo Kawaguchi, \emph{The gonality and the Clifford index of curves on a toric surface}, preprint
 \bibitem{Khovanskii} Askold G.\ Khovanskii, \emph{Newton polyhedra and toroidal varieties},
Functional Analysis and Its Applications \textbf{11}(4), pp.\ 289-296 (1977)
 \bibitem{Koelman} Robert J.\ Koelman, \emph{The number of moduli of families of curves on toric surfaces},
 Ph.D. thesis, Katholieke Universiteit Nijmegen (1991)
 \bibitem{Koelman2} Robert J.\ Koelman, \emph{A criterion for the ideal of a projectively embedded toric surface to be generated
  by quadrics}, Beitr\"age zur Algebra und Geometrie \textbf{34}, pp.\ 57-62 (1993)
\bibitem{kollar} J\'anos Koll\'ar, Shigefumi Mori, \emph{Birational geometry of algebraic varieties}, Cambridge Tracts in Mathematics \textbf{134}, Cambridge University Press (1998)
\bibitem{LelliChiesa} Margherita Lelli-Chiesa, \emph{Green's conjecture for curves on surfaces with an anticanonical pencil}, Mathematische Zeitschrift \textbf{275}, pp.\ 899-910 (2013)
 \bibitem{LubbesSchicho} Niels Lubbes, Josef Schicho, \emph{Lattice polygons and families of curves on rational surfaces}, Journal of Algebraic
 Combinatorics \textbf{34}, pp.\ 213-236 (2012)
 \bibitem{Matsumoto} Ryutaroh Matsumoto, \emph{The $C_{ab}$ curve}, note available at \verb"http://www.rmatsumoto.org/cab.pdf"
 \bibitem{Miura} Shinji Miura, \emph{Error-correcting codes based on algebraic geometry}, Ph.D.\ thesis, University of Tokyo (1997)
 \bibitem{Namba} Makoto Namba, \emph{Families of meromorphic functions on compact Riemann surfaces},
    Lecture Notes in Mathematics \textbf{767}, Springer (1979)
\bibitem{reid} Miles Reid, \emph{Nonnormal del Pezzo surfaces}, Publications of the Research Institute for Mathematical Sciences \textbf{30}(5), pp.\ 695-727 (1994)
 \bibitem{saintdonat} Bernard Saint-Donat, \emph{On Petri's analysis of the linear system of quadrics through a canonical curve}, Mathematische 
  Annalen \textbf{206}, pp.\ 157-175 (1973)
\bibitem{Schenck} Hal Schenck, \emph{Lattice polygons and Green's theorem}, Proceedings of the American Mathematical Society \textbf{132}(12), pp.\ 3509-3512 (2004)
 \bibitem{Schreyer} Frank-Olaf Schreyer, \emph{Syzygies of canonical curves and special linear series}, Mathematische Annalen \textbf{275}(1), pp.\ 105-137 (1986)
 \bibitem{Serrano} Fernando Serrano, \emph{Extensions of morphisms defined on a divisor}, Mathematische Annalen \textbf{277}(3), pp. 395-413 (1987)
\end{thebibliography}
\end{document}